\documentclass[preprint,10pt,a4paper]{elsarticle}

\usepackage[margin=1in]{geometry}

\usepackage[utf8]{inputenc}

\usepackage{amsmath}
\usepackage{amsthm}
\usepackage{amssymb}
\usepackage{amsfonts}
\usepackage{a4wide}
\usepackage{xcolor}
\usepackage{subcaption}
\usepackage[colorinlistoftodos]{todonotes}
\usepackage{bm}
\usepackage[T1]{fontenc}
\usepackage{graphicx}
\usepackage{colortbl,dcolumn}

\journal{t.b.d.}

\begin{document}

\begin{frontmatter}



\title{Nonlinear model order reduction for problems with microstructure \\ using mesh informed neural networks}


\author[1]{Piermario Vitullo}
\author[2]{Alessio Colombo}
\author[1]{Nicola Rares Franco}
\author[1]{Andrea Manzoni}
\author[1]{Paolo Zunino}

\address[1]{MOX, Department of Mathematics, Politecnico di Milano, P.zza Leonardo da Vinci 32, 20133, Milan, Italy}

\address[2]{Department of Civil and Environmental Engineering, Politecnico di Milano, P.zza Leonardo da Vinci 32, 20133, Milan, Italy}

\begin{abstract}
Many applications in computational physics involve approximating problems with microstructure, characterized by multiple spatial scales in their data. However, these numerical solutions are often computationally expensive due to the need to capture fine details at small scales. As a result, simulating such phenomena becomes unaffordable for many-query applications, such as parametrized systems with multiple scale-dependent features. Traditional projection-based reduced order models (ROMs) fail to resolve these issues, even for second-order elliptic PDEs commonly found in engineering applications. To address this, we propose an alternative nonintrusive strategy to build a ROM, that combines classical proper orthogonal decomposition (POD) with a suitable neural network (NN) model to account for the small scales. Specifically, we employ sparse mesh-informed neural networks (MINNs), which handle both spatial dependencies in the solutions and model parameters simultaneously. We evaluate the performance of this strategy on benchmark problems and then apply it to approximate a real-life problem involving the impact of microcirculation in transport phenomena through the tissue microenvironment.
\end{abstract}


% \begin{highlights}
% \end{highlights}


\end{frontmatter}


%%%%%%%%%%%%%%%%%%%%%%%%%%%%%%%%%%%%%%%%%%%%%%%%%%%%%%%%%%%%%%%%%%%%%%
%User defined commnands

\def\eg{\textit{e.g.}\ }
\def\ie{\textit{i.e.}\ }
\def\argmin{\mathrm{argmin}}
\def\esssup{\mathrm{ess\,sup}}
\def\I{\boldsymbol{\mathrm{I}}}
\def\Tau{\mathcal{T}}
\def\para{\boldsymbol{\mu}}
\def\parmi{\boldsymbol{\mu}_{m}}
\def\parma{\boldsymbol{\mu}_{M}}
\def\micro{m}
\def\macro{M}
\def\parph{\boldsymbol{\mu}_{ph}}
\def\parg{\boldsymbol{\mu}_{g}}
\def\domega{\partial \Omega}
\def\uph{u_{\para,h}}
\def\uuphi{\mathtt{\bf u}_{\para^{(i)},h}}
\def\uuh{\mathtt{\bf u}_{h}}
\def\uuph{\mathtt{\bf u}_{\para,h}}
\def\uphrb{u_{\para,h}^{rb}}
\def\nuphrb{\widetilde{u}_{\para,h}^{rb}}
\def\nuuphrb{\widetilde{\bf u}_{\para,h}^{rb}}
\def\muphrb{\widehat{u}_{\para,h}^{rb}}
\def\muuphrb{\widehat{\mathbf{u}}_{\para,h}^{rb}}
\def\muuphrbj{\widehat{\mathbf{u}}_{\para,h}^{rb,(j)}}
\def\uuphrb{\mathtt{\bf u}_{\para,h}^{rb}}
\def\urb{\mathtt{u}_{\para,rb}}
\def\nurb{\widetilde{\mathtt{u}}_{\para,rb}}
\def\nuu{\widetilde{\mathtt{u}}_{\para}}
\def\nuurb{\widetilde{\mathtt{\bf u}}_{\para,rb}}
\def\nuuu{\widetilde{\mathtt{\bf u}}_{\para,c}}
\def\muurb{\widehat{\mathtt{\bf u}}_{\para,rb}}
\def\uurb{\mathtt{\bf u}_{\para,rb}}
\def\npmicro{n_{\micro}}
\def\npmacro{n_{\macro}}
\def\nrb{n_{rb}}
\def\nnrb{\widetilde{n}_{rb}}
\def\npsi{\widetilde{\psi}}
\def\nft{\widetilde{\alpha}}
\def\RR{\mathbb{R}}
\def\VV{\mathbb{V}}
\def\WW{\mathbb{W}}
\def\SS{\mathbb{S}}
\def\ann{\mathcal{N}}
\def\lnn{\mathcal{L}}
\def\minn{\mathcal{M}}
\def\hyper{\boldsymbol{\theta}}
\def\hypermi{\boldsymbol{\theta}_{\micro}}
\def\hyperd{\boldsymbol{\theta}_{d}}
\def\hyperi{\boldsymbol{\theta}_{\eta}}
\def\hyperma{\boldsymbol{\theta}_{\macro}}
\def\wwmic{\boldsymbol{\widetilde{\theta}}_{g}}
\def\uu{\mathbf{v}}
\def\inlets{\boldsymbol{\eta}}
\def\dist{\mathbf{d}}
\newcommand{\bct}{{C}_t}
\newcommand{\bp}{{p}}

\theoremstyle{definition}
\newtheorem{definition}{Definition}[section]

%%%%%%%%%%%%%%%%%%%%%%%%%%%%%%%%%%%%%%%%%%%%%%%%%%%%%%%%%%%%%%%%%%%%%%
\section{Introduction}
\label{sec:intro}

The repeated solution of differential problems required to describe, forecast, and control the behavior of a system in multiple virtual scenarios is a computationally extensive task if relying on classical, high-fidelity full order models (FOMs) such as, e.g., the finite element method, when very fine spatial grids and/or time discretizations are employed to capture the detailed behavior of the solution. The presence of multiple (spatial and/or temporal) scales affecting problem's data and its input parameters -- such as, e.g., material properties or source terms included in the physical model -- makes this task even more involved, and is a common issue whenever dealing with, e.g., biological models \cite{El-Bouri201540,Vidotto20191076,Possenti2019}, structural mechanics \cite{Meier2018124,Fumagalli2013454} as well as environmental flows \cite{Tschisgale2021,Hagmeyer2022}, just to make a few examples.
Reduced order modeling techniques provide nowadays a wide set of numerical strategies to solve parametrized differential problems in a rapid and reliable way. 
For instance, in the case of problems where the microstructure is characterized by slender fibers immersed into a continuum, reduced order models based on dimensional reduction of the fibers have been successfully adopted, in the framework of mixed dimensional problem formulations, see, e.g.,  \cite{D'Angelo20081481,Kuchta2021558,Mori2019887}.
However, given the increasing complexity of the problems that need to be addressed, mixed dimensional formulations are no longer sufficient; indeed, they actually represent the starting point for a second level of model reduction that is addressed by our work.

Physics-based ROMs such as, e.g., the reduced basis (RB) method \cite{Hestaven1,QuartManzNegri}, provide a mathematically rigorous framework to build ROMs involving a linear trial subspace, or reduced basis -- obtained, e.g., through proper orthogonal decomposition (POD) on a set of FOM snapshots -- to approximate the solution manifold, and a (Petrov-)Galerkin projection to derive the expression of the ROM. Provided the problem operators depend affinely on the input parameters, a suitable offline-online splitting ensures a very fast assembling and solution of the ROM during the online testing stage, once the linear subspace and the projected reduced operators have been computed and stored during the offline training stage. 

The presence of involved, spatially-dependent, input parameters representing, e.g., diffusivity fields or distributed force fields within the domain where the problem is set, usually makes the offline-online splitting not so straightforward because of their high dimension and (potentially, highly) nonlinear nature of the parameter-to-solution map. To be consistent with the usual formulation of ROMs, we will refer, throughout the paper, to space-varying fields representing some problem inputs as to {\em input parameters}, despite these latter usually denote vectors of quantities. These features might impact on the reducibility of problems with microstructure with {\em linear} ROMs at least in two ways: {\em (i)} linear trial subspaces might have large dimension, thus not really reducing the complexity of the problem; %compared to the intrinsic dimension of the solution manifold being approximated;(IN REALTA', qui anche la dimensione intrinseca è grande, perché lo spazio parametrico ha -a priori- dimensione elevata)
{\em (ii)} classical hyper-reduction techniques, such as the (discrete) empirical interpolation method (DEIM) \cite{barrault2004empirical,chaturantabut2010nonlinear,negri2015efficient,bonomi2017reduced} or the energy-conserving sampling and weighting (ECSW) method \cite{farhat2015structure,farhat20205}, usually employed to speed up the ROM assembling, can also suffer from severe computational costs, and result in intrusive procedures.  

To overcome these drawbacks, alternative data-driven methods can be used for the approximation of the RB coefficients without resorting to (Galerkin) projection. In these cases, the FOM solution is projected onto the RB space and the combination coefficients are approximated by means of a surrogate model, exploiting, e.g., radial basis function (RBF) approximation \cite{audouze2013nonintrusive}, polynomial chaos expansions \cite{sun2021non}, artificial neural networks (NNs) \cite{hesthaven2018non, gao2021non, salvador20211}, or Gaussian process regression (GPR) \cite{guo2018reduced, guo2019, kast2020non, zhang2019model}. The high-fidelity solver is thus used only offline to generate the data required to build the reduced basis, and then to train the surrogate model. Non-intrusive POD-based approaches using RBFs to interpolate the ROM coefficients have been proposed in, e.g., \cite{audouze2013nonintrusive, walton2013reduced, xiao2017parameterized}. In a seminal contribution, Hesthaven and Ubbiali \cite{hesthaven2018non} have instead employed NNs to build a regression model to compute the coefficients of a POD-based ROM, in the case of steady PDEs; a further extension to time-dependent nonlinear problems, i.e., unsteady flows, has been addressed in \cite{wang2019non}; see also, e.g., \cite{swischuk2019projection,yu2019non} for the use of alternative machine learning strategies to approximate the POD coefficients. We will refer hereon to this class of approaches as to {\em POD-NN} methods. POD with GP regression have been used to build ROMs by Guo et al. in \cite{guo2018reduced} to address steady nonlinear structural analysis, as well as for time-dependent problems  \cite{guo2019}. A detailed comparison among non-intrusive ROMs employing RBF, ANN, and GP regressions can  be found in, e.g., \cite{berzicnvs2020standardized}; see instead \cite{guo2022multi,conti2022multi} for an alternative use of NNs to perform regression in the context of multi-fidelity methods, capable of leveraging models with different fidelities. These latter may involve data-driven projection-based \cite{peherstorfer2016data, guo2022bayesian} ROMs or more recently developed deep learning-based ROMs \cite{fresca2021comprehensive,fresca2021pod, FrancoROMPDE, botteghi2022deep}. In these latter cases, POD has been replaced by (e.g., convolutional) autoencoders to enhance dimensionality reduction of the solution manifold, relying on (deep) feedforward NNs to learn the reduced dynamics onto the reduced trial manifold. 

Despite the advantages they provide compared with dense architectures -- in terms of costs and size of the optimization problem to be solved during training -- convolutional NNs cannot handle general geometries and they might become inappropriate as soon as the domain where the problem is set is not an hypercube, although some preliminary attempts to generalize CNN in this direction have recently appeared \cite{gao2021phygeonet}. In the case of problems with microstructure, this issue usually arises when attempting at reducing the dimensionality of spatially distributed parameters using DNNs rather than POD as, e.g., the DEIM would do -- provided a linear subspace built through POD is still employed to reduce the solution manifold. For this reason, in this paper we rely on Mesh-Informed Neural Networks (MINNs), a class of architectures recently introduced in \cite{FrancoMINN} and specifically tailored to handle mesh based functional data. The driving idea behind MINNs is to embed hidden layers into discrete functional spaces of increasing complexity, obtained through a sequence of meshes defined over the underlying spatial domain. This approach leads to a natural pruning strategy, whose purpose is to inform the model with the geometrical knowledge coming from the domain. As shown through several examples in \cite{FrancoMINN}, MINNs offer reduced training times and better generalization capabilities, making them a competitive alternative to other operator learning approaches, such as DeepONets \cite{lu2021learning} and Fourier Neural Operators \cite{li2020fourier}. Our purpose is to employ MINNs to enable the design of sparse architectures aiming at feature extraction from space-varying parameters that define the problem's microstructure. 

In this paper we propose a new strategy to tackle parametrized problems with microstructure, combining a POD-NN method to build a reduced order model in a nonintrusive way and MINNs to build a closure model capable of integrating in the resulting approximation the information coming from mesh based functional data, in order to avoid to deal with a very large number of POD modes in presence of complex microstructures. Such a problem arises, e.g., when describing  oxygen transfer in the microcirculation by including blood flow and hematocrit transport coupled
with the interstitial flow, oxygen transport in the blood and the tissue, described by the vascular-tissue exchange. The presence of microvasculature, described in terms of a (varying) graph-like structure within the domain, requires the use of a NN-based strategy to reduce the dimensionality of such data, thus calling into play MINNs, yielding a strategy we refer to as a POD-MINN method. 
To account for the neglected scales at the POD level, and correct the POD-MINN approximation to enhance its accuracy without further increasing its dimension, we equip the POD-MINN method with a closure model, ultimately yielding a strategy we refer to as a POD-MINN+ approximation. \\

The structure of the paper is as follows. In Sect.~\ref{sec:setup} we formulate the class of problems with microstructure we focus on in this paper, introducing their high-fidelity approximation by the finite element method, recalling how classical projection-based ROMs are formulated, and showing their main limitations. In Sect.~\ref{sec:nn} we describe how to take advantage of mesh-informed neural networks to tackle problems with microstructure, while in Sect.~\ref{sec:pod-nn} we address the POD-MINN and the POD-MINN+ methods, showing results obtained in a series of simple numerical test cases. Finally, in Sect.~\ref{sec:apps} we consider the application of the proposed strategy to an oxygen transfer problem taking place in the microcirculation, and draw some conclusions in Sect.~\ref{sec:conclusions}.

%%%%%%%%%%%%%%%%%%%%%%%%%%%%%%%%%%%%%%%%%%%%%%%%%%%%%%%%%%%%%%%%%%%%%%
\section{Formulation and approximation of PDEs with embedded microstructure}
\label{sec:setup}

\subsection{Problem setup}

%We describe here in general terms a problem 
We start by describing, in general terms, the essential properties of a problem governed by a parametrized PDE affected by a microstructure. With the term \emph{microstructure} we refer to some features, primarily of the forcing terms of the equations, that induce the coexistence of multiple characteristic scales in the solution. Despite this is a particular case of \emph{multiscale} problem, developing ROMs for multiscale problems in the spirit of upscaling and/or numerical homogenization is not the scope of this work -- this aspect has been addressed, at least under some simplifying assumptions, by several authors in the framework of reduced basis methods, see, e.g., \cite{boyaval2008reduced,abdulle2012reduced,abdulle2014offline}. Conversely, here we aim to develop a ROM that fully captures all the scales of the solution.

We restrict to steady problems governed by second order elliptic equations and we assume that the microstructure may influence the parameters of the operator, the forcing terms and the boundary conditions, whereas the domain is fixed and do not depend on parameters. 
Under these assumptions, in this paper we address parametrized PDEs of the form 
\begin{equation}\label{eq:problem}
    \begin{cases}
    L_{\para} u_{\para} = f_{\para} & \text{in} \ \Omega,
    \\
    B_{\para} u_{\para}= g_{\para} & \text{on} \ \domega,
\end{cases}
\end{equation}
where the solution $u_{\para}$ depends on the parameter vector $\para = (\parma, \parmi) \in \mathcal{P}=\mathcal{P}_{\macro}\times\mathcal{P}_{\micro}$, with $\mathrm{dim}(\mathcal{P}_{\macro})=\npmacro$ and 
$\mathrm{dim}(\mathcal{P}_{\micro})=\npmicro$. In other words, $u_{\para}=u(\para)$. Both notations will be used in what follows, with a preference to the most compact one when the context allows it.
%
 %The inputs of the problem may be discrete or spatially distributed functions. 
%Since the beginning
%
By the distinction between \emph{macroscale parameters} -- denoted by $\parma$ -- and \emph{microscale parameters} -- denoted instead by $\parmi$ -- we highlight that problem parameters satisfy a {\em scale separation property}. 
%We subdivide accordingly the parameter space in two subspaces, corresponding to the \emph{micro} and \emph{macro}-scales respectively. Precisely, we assume that 
%For the problem addressed asapplication of the proposed methodology,
For example, in the biophysical application discussed at the end of this work, the physical parameters of the operator are affected by the macroscale parameter $\parma$, namely $L_{\para}=L_{\parma}$, while the geometry of the microstructure determines the forcing terms of the problem, that is, $f_{\para}=f_{\parmi}$ and $g_{\para}=g_{\parmi}$. 
As a consequence of these assumptions, we focus on the particular case of the abstract problem \eqref{eq:problem}, that can be rewritten as follows,
\begin{equation}\label{eq:micro_problem}
    \begin{cases}
    L_{\parma} u_{\para} = f_{\parmi} & \text{in} \ \Omega,
    \\
    B_{\parmi} u_{\para}= g_{\parmi} & \text{on} \ \domega.
\end{cases}
\end{equation}

In this particular context, we implicitly assume that the physical parameters belong to a low-dimensional space, whereas the geometry of the microstructure features high dimensionality -- that is, $\npmicro\gg \npmacro$. Moreover, for the sake of presentation, in this section we assume that the operator $L_{\parma}$ appearing in problem \eqref{eq:micro_problem} is linear, although our methodology will also be applied to a nonlinear oxygen transfer problem, whose formulation is briefly described in Sect.~\ref{sec:micromicro}. % in the applications the proposed methodologies, to the nonlinear case. Nevertheless, we will consider problems where the dependence of the operators and forcing terms from the parameters is, in general, non-affine.

%%%%%%%%%%%%%%%%%%%%%%%%%%%%%%%%%%%%%%%%%%%%%%%%%%%%%%%%%%%%%%%%%%%%%%
\subsection{An example of problem with microstructure inspired to microcirculation}\label{sec:micromicro}

We also address in this paper a mathematical model for oxygen transfer in the microcirculation on the basis of a comprehensive model described in \cite{Possenti2019,Possenti20213356} that includes blood flow and hematocrit transport coupled with the interstitial flow, oxygen transport in the blood and the tissue, described by the vascular-tissue exchange. Indeed, our proposed methodology can be applied to any field described by the general microcirculation model, such as fluid pressure, velocity, oxygen concentration (or partial pressure), despite we only focus, in this paper, on the oxygen transport. The non-intrusive character of the POD-MINN (and POD-MINN+) method indeed allows us to approximate a single field, despite the coupled nature of the general model. Note that if we relied on a classical projection-based ROM we had to approximate all the variables simultaneously, thus requiring to face a much higher degree of complexity when constructing the ROM. 
 
The general model describes the flow in two different domains,  the tissue domain ($\Omega \subset \RR^3$ with $\dim(\Omega)=3$), where the  unknowns are the fluid pressure $p_t$, the fluid velocity $\uu_t$, and the oxygen concentration $C_t$, and the vascular domain ($\Lambda \subset \RR^3$ with $\dim(\Lambda)=1$), which is a metric graph describing a network of connected one-dimensional channels, where the unknowns are the blood pressure $p_v$, the blood velocity $\uu_v$, and the vascular oxygen concentration $C_v$. The model for the oxygen transport employs the velocity fields $\uu_v$ and $\uu_t$ to describe blood flow in the vascular network and the plasma flow in the tissue. The governing equations of the oxygen transfer model read as follows:
\begin{equation}
\label{eq:oxy}
\left\{
\begin{array}{ll} 
\displaystyle
\nabla \cdot \left(-D_t \nabla C_t + \uu_t~ C_t\right) 
+ V_{max} ~ \frac{C_{t}}{C_{t} + \alpha_t ~p_{m_{50}}} = \phi_{O_2} \, \delta_\Lambda \quad &\textrm{on $\Omega$}
\\[6pt]       
\displaystyle
\pi R^2 \frac{\partial}{\partial s} \left(
- D_v \frac{\partial C_v}{\partial s} + {v_v} ~ C_v 
+ {v_v} ~ k_1~H~  \frac{C_v^\gamma}{C_v^\gamma + k_2}\right)  
= -  \phi_{O_2} \quad &\textrm{on $\Lambda$}
\\[6pt]
\displaystyle
\phi_{O_2} = 2 \pi R~P_{O_2} (C_v - \bct) + (1-\sigma_{O_2})~\left(\frac{C_v + \bct}{2}\right)~\phi_v \quad &\textrm{on $\Lambda$}
\\[10pt]
\displaystyle
\phi_v=2 \pi R L_p\big((p_v - \bp_t) - \sigma(\pi_v-\pi_t)\big)
\\[6pt]
\displaystyle
C_v = C_{in}  \quad &\textrm{on}~\partial\Lambda_{\text{in}}
\\[6pt]
\displaystyle
- D_v \frac{\partial C_v}{\partial s} = 0  \quad &\textrm{on}~\partial\Lambda_{\text{out}}
\\[6pt]
\displaystyle
- D_t \nabla C_t\cdot\mathbf{n} = \beta_{O_2}~(C_t - C_{0,t}) \quad &\textrm{on}~\partial\Omega.
\end{array}
\right.
\end{equation}
In particular, the first equation governs the oxygen in the tissue, the second describes how the oxygen is transported by the blood stream, and the third defines the oxygen transfer between the two domains $\Omega$ and $\Lambda$. In particular, the model for the flux $\phi_{O_2}$ is obtained assuming that  the vascular wall acts as a semipermeable membrane. This model is complemented with a set of boundary conditions reported in the last three equations: at the vascular inlets $\partial\Lambda_{\text{in}}$ we prescribe the oxygen concentration; at the vascular endpoints $\partial\Lambda_{\text{out}}$ null diffusive flux is enforced; and for the boundary of the tissue domain $\partial\Omega$ we simulate the presence of an adjacent tissue domain with a boundary conductivity $\beta_{O_2}$ and a far-field concentration $C_{0,t}$.
The symbols $D_t,D_v,V_{max},\alpha_t,p_{m_{50}},k_1,k_2,C_v^\gamma,P_{O_2},\sigma_{O_2},L_p,\sigma,\pi_v,\pi_t$ are constants independent of the solution of the model. For a detailed description of the physical meaning of these quantities see, e.g.,  \cite{Possenti20213356}.

Comparing model \eqref{eq:micro_problem} with \eqref{eq:oxy}, we observe that the operator $L_{\parma}$ consists of the left hand side of the first equation, where the macroscale parameters are the physical parameters of the operator $\nabla \cdot \left(-D_t \nabla C_t + \uu_t~ C_t\right) + V_{max}C_t/(C_{t} + \alpha_t ~p_{m_{50}})$, such as for example $V_{max}$ that modulates the oxygen consumption by the cells in the tissue domain. Because of the last term, such operator is nonlinear. The solution $u_{\para}$ corresponds to the concentration $C_t$, and the flux $\phi_{O_2} \delta_\Lambda$ plays the role of the forcing term $f_{\parmi}$. As $\delta_\Lambda$ denotes a Dirac  $\delta$-function distributed on $\Lambda$, the term $\phi_{O_2} \delta_\Lambda$ represents an incoming mass flux for the oxygen concentration $C_t$, supported on the vascular network $\Lambda$. In other words, it is a concentrated source defined in $\Omega$, acting as a microscale forcing term. According to the third equation, the intensity of the oxygen flux depends on both concentrations, namely $\phi_{O_2}=\phi_{O_2}(C_v,C_t)$. In turn, the vascular oxygen concentration $C_v$ is governed by the second equation and by the inflow boundary conditions $C_v = C_{in}$ on $\partial\Lambda_{IN}$. This condition plays the role of $B_{\parmi} u_{\para}= g_{\parmi}$ on $\domega$ in problem \eqref{eq:micro_problem}, and as it depends on the inlet points of the vascular network on the domain boundary, it is classified as a microscale model. In conclusion, we can establish the following connections between problems \eqref{eq:micro_problem} and \eqref{eq:oxy}:
\begin{align*}
    u_{\para} &\approx C_t,
    \\
    L_{\parma}u_{\para} &\approx 
    \nabla \cdot \left(-D_t \nabla C_t + \uu_t~ C_t\right) 
    + V_{max} ~ \frac{C_{t}}{C_{t} + \alpha_t ~p_{m_{50}}},
    \\
    f_{\parmi} &\approx \phi_{O_2}(C_v,C_t)\delta_\Lambda,
    \\
    B_{\parmi} u_{\para} &\approx C_v|_{\partial\Lambda_{IN}},
    \\
    g_{\parmi} &\approx C_{in}|_{\partial\Lambda_{IN}}.
\end{align*}

We remark that problem \eqref{eq:oxy} is just an example among many other relevant problems characterized by a microstructure, such as flows through perforated or fractured domains, or the mechanics of fiber-reinforced materials, just to mention a few. Our proposed methodology does not depend on the specific problem and can be applied with suitable adaptations to all these problems thanks to its non-intrusive nature.

%%%%%%%%%%%%%%%%%%%%%%%%%%%%%%%%%%%%%%%%%%%%%%%%%%%%%%%%%%%%%%%%%%%%%%
\subsection{Numerical approximation by the finite element method}

We consider the finite element method (FEM) as the high-fidelity FOM for the problem described in the previous section, as well as for the other test cases we propose. The central hypothesis underlying this work is that the FEM approximation of the microscale problem resolving all the scales of the microstructure is computationally too expensive for real-life applications. For the sake of simplicity, we present our methodology referring to the general, abstract problem \eqref{eq:problem}. %We name such an approach the full order model and the main motivation of this work is to propose a reduced order model that significantly decreases the computational cost of the former.
Before introducing its finite element discretization, we address the variational formulation of \eqref{eq:problem} that, given a suitable Hilbert space $(V, \| \cdot \|)$ depending on the boundary conditions imposed on the problem at hand, reads as follows: for any $\para \in \mathcal{P} \subset \RR^{n_p}$, find $u_{\para} \in V$ such that
\begin{equation}\label{eq:problem_weak}
a_{\para}(u_{\para}, v) = F_{\para}(v), \quad \forall v \in V,
\end{equation}
being $a_{\para}: V \times V \mapsto \mathbb{R}$ and $F_{\para}: V \mapsto \mathbb{R}$ two parameter-dependent operators. %Problem \eqref{eq:problem} (or \eqref{eq:problem_weak}) therefore defines
Depending on the formulation of choice, Problems \eqref{eq:problem} and \eqref{eq:problem_weak} naturally define a \textit{parameter-to-solution map}, that is,
a map that assigns to each input parameter $\para \in \mathcal{P}$ the corresponding solution $u_{\para}=u(\para)$, \ie, 
\begin{equation*}
    u_{\para}: \, \mathcal{P} \,  \mapsto \, V;
    \quad
    \para         \,  \mapsto     \, u(\para).
\end{equation*}
This also allows us to define the \textit{solution manifold} $\mathcal{S} = \{u_{\para}=u(\para) : \para \in \mathcal{P}\}$. 

The high-fidelity FOM consists of the Galerkin projection of problem \eqref{eq:problem_weak} onto a Finite Element (FE) space $V_h$ of dimension $N_h=\mathrm{dim}(V_h)$, suitably chosen depending on the characteristics of the problem at hand. Assuming for simplicity a fully conformal approximation, given $V_h \subset V$, we aim to find $\uph \in V_h$ such that
\begin{equation}\label{eq:problem_weak_h}
a_{\para}(\uph, v_h) = F_{\para}(v_h) \qquad \forall v_h \in V_h.
\end{equation}
As pointed out for the continuous case, problem \eqref{eq:problem_weak_h} defines a mapping $\mathcal{P} \, \mapsto \, V_h$ that identifies the discrete solution manifold $\mathcal{S}_h = \{\uph: \para \in \mathcal{P}\}$.
From the discrete standpoint, problem \eqref{eq:problem_weak_h} is equivalent to a (large) system of algebraic equations of the form
\begin{equation*}
    A_{\para,h}\uuph = \mathbf{F}_{\para,h}
\end{equation*}
where $\uuph \in \mathbb{R}^{N_h}$ is the vector of degrees of freedom of the FE approximation. The need to, e.g., sample adequately the statistical variability of the microstructure, and assess its impact on the problem solution, would ultimately require to repeatedly query such FOM, whence the need of a rapid and reliable ROM. %. Recuced order modeling \cite{QuartManzNegri} sets out to provide a solution to this problem by replacing the full order model with one that is computationally cheaper.

%%%%%%%%%%%%%%%%%%%%%%%%%%%%%%%%%%%%%%%%%%%%%%%%%%%%%%%%%%%%%%%%%%%%%%
\subsection{Projection-based model order reduction}

For second order elliptic PDEs, projection-based ROMs represent a consolidated approach that has been successful in many areas of application \cite{benner2015survey}. The main idea is to generate the ROM by projecting the FOM onto a low-dimensional linear subspace of $V_h$. Precisely, the \textit{reduced basis} (RB) method aims at building a subspace $V_{rb} \subset V_h$ through the linear combination of a set of $\nrb \ll N_h$ basis functions, being $\nrb$ independent of $h$ and $N_h$:
\begin{equation*}
    V_{rb} = \mathrm{span}\{\psi_1, \ldots, \psi_{\nrb}\}.
\end{equation*}
%The RB approximation $\uphrb \in V_{rb}$ is then defined as the linear combination of the reduced basis
%\begin{equation}\label{eq:lin_reduced_basis}
%    \uphrb = \sum_{i=1}^{\nrb} \urb^{(i)}\psi_i(x) 
%\end{equation}
%where $\uurb=[\urb^{(1)},\ldots \urb^{(i)}, \ldots \urb^{(\nrb)}]^T$ are the \emph{reduced coefficients}, and $\uphrb$ is determined by the Galerkin projection of $\uph$ on $V_{rb}$; precisely,
%\begin{equation*}
%    a_{\para} (\uphrb,v_{rb})=F_{\para}(v_{rb}), \qquad \forall v_{rb} \in V_{rb}.
%\end{equation*}
%corresponding the following discrete algebraic system,
%\begin{equation}\label{eq:system_ROM}
%    A_{\para,rb} \uurb = \mathbf{F}_{\para,rb}.
%\end{equation}
The RB approximation $\uphrb \in V_{rb}$ is then defined as the solution to the Galerkin projection of the FOM onto the RB space, that is,
\begin{equation}
    \label{eq:RB-pde}
    a_{\para} (\uphrb,v_{rb})=F_{\para}(v_{rb}), \qquad \forall v_{rb} \in V_{rb}.
\end{equation}
Equivalently, by expanding the RB solution over the reduced basis, the latter can be expressed as
\begin{equation}\label{eq:lin_reduced_basis}
    \uphrb = \sum_{i=1}^{\nrb} \urb^{(i)}\psi_i(x) 
\end{equation}
where $\uurb=[\urb^{(1)},\ldots,\urb^{(\nrb)}]^T$ are the \emph{reduced coefficients}, which, in practice, are obtained by solving the algebraic counterpart of \eqref{eq:RB-pde}, namely
\begin{equation}\label{eq:system_ROM}
    A_{\para,rb} \uurb = \mathbf{F}_{\para,rb}.
\end{equation}

Overall, equations \eqref{eq:RB-pde}, \eqref{eq:lin_reduced_basis} and \eqref{eq:system_ROM} illustrate different ways to represent the ROM approximation of the FOM solution. Precisely, $\uurb \in \RR^{\nrb}$ are the reduced coefficients that identify the ROM solution while $\uuphrb \in \RR^{N_h}$ is the representation of the same solution in the high-fidelity (FEM) space. We note that, although they have the same dimension, $\uuph$ and $\uuphrb$ do not coincide.

We stress that this projection determines an approximation of the FOM solution $\uph$ onto a \emph{linear subspace} of the discrete solution manifold $\mathcal{S}_h$. As a consequence the corresponding error behaves -- at least, theoretically -- following the decay of the Kolmogorov $\nrb$-width \cite{QuartManzNegri}. Indeed, a fast-decaying Kolmogorov $\nrb$-width reflects the approximability of the solution manifold by finite-dimensional linear spaces.

From an algebraic standpoint, the reduced basis is encoded in the matrix $\mathbb{V} = [\boldsymbol {\psi}_1 \ | \ \ldots \ | \ \boldsymbol {\psi}_{\nrb} ] \in \mathbb{R}^{N_h \times \nrb}$, where $\boldsymbol {\psi}_i \in \mathbb{R}^{N_h}$ are the FOM degrees of freedom of the $i$-th basis $\psi_i$, so that \eqref{eq:lin_reduced_basis} can be equivalently rewritten as $\uuphrb= \mathbb{V} \uurb$ and the RB problem can be assembled projecting the FOM matrix and right hand side onto the reduced space,
\begin{equation}\label{eq:projection_FOM}
   A_{\para,rb} = \mathbb{V}^T A_{\para,h} \mathbb{V},
   \quad
   \mathbf{F}_{\para,rb} = \mathbb{V}^T \mathbf{F}_{\para,h}.
\end{equation}
For a special category of problems, characterized by an affine parameter dependence, the algorithm presented to build and solve \eqref{eq:system_ROM} can be  split into a parameter-independent \emph{offline} and a parameter-dependent \emph{online} stage \cite{QuartManzNegri}. However, if the affine parameter dependence does not hold, suitable hyper reduction techniques must be called into play to restore it, at least in an approximate way; nevertheless, this is not a focus of this work. 

Regarding the construction of the subspace $V_{rb}$, a classical choice is to rely on POD, starting from  % with good (possibly optimal) approximation properties. To this aim there is a vast literature embracing two main families of approaches, the \emph{proper orthogonal decomposition} and the \emph{greedy methods}. We opt for the former, although the main methodologies proposed in this work are independent of this particular choice.
%
%To define a reduced basis using Proper Orthogonal Decomposition (POD), we consider 
a collection of FOM solutions, named \emph{snapshots}, for suitably chosen parameters $\para^{(i)}$, such that the span of these functions, $\mathrm{span}\{u_h(\para^{(i)})\}_{i=1}^N$ is a reasonable approximation of the discrete solution manifold $\mathcal{S}_h$.
By collecting the snapshots FOM degrees of freedom as the columns of a \emph{data matrix} %$\mathbb{S}$, that is 
\begin{equation*}
\mathbb{S}=[\mathtt{\bf u}_h(\para^{(1)})|\ldots | \mathtt{\bf u}_h(\para^{(i)})| \ldots|\mathtt{\bf u}_h(\para^{(N)})],
\end{equation*}
the singular value decomposition of $\mathbb{S}$,
\[
\mathbb{S}=\mathbb{W} \left[\begin{array}{cc}
\widetilde{\mathbb{D}} &  \mathbf{0} \\
\mathbf{0} & \mathbf{0}
\end{array} \right] \mathbb{Z}^T
\]
ensures the existence of three matrices $\mathbb{W} \in \mathbb{R}^{N_h \times N_h}$, $\widetilde{\mathbb{D}}=\mathrm{diag}[\sigma_1,\ldots,\sigma_R]$, being $\sigma_1\geq\ldots\geq\sigma_i\geq\ldots\geq\sigma_R$ the singular values of $\mathbb{S}$ and $R$ is its rank,
and $\mathbb{Z} \in \mathbb{R}^{N \times N}$. The columns of $\mathbb{W}$ and $\mathbb{Z}$ are named the left and right singular vectors, respectively. The POD basis is then defined as the set of the first $\nrb$ left singular vectors of $\mathbb{S}$; their collection provides the matrix $\widetilde{\mathbb{W}}(\nrb)$ that represents the \emph{best} $\nrb$-rank approximation of $\mathbb{S}$ with respect to the following error \cite{QuartManzNegri},
\begin{equation*}
    E(\widetilde{\mathbb{W}}(\nrb);\mathbb{S})=\min\limits_{\mathbb{M} \in \mathbb{R}^{N_h \times N_h}} \sum_{j=1}^N \|\mathbb{S}(:,j) - \mathbb{M} \mathbb{M}^T  \mathbb{S}(:,j) \|_{\mathbb{R}^{N_h}}^2.
     %E(\widetilde{\mathbb{W}}(\nrb);\mathbb{S})=\min\limits_{\mathbb{M} \in \mathbb{R}^{N_h \times N_h}} \sum_{j=1}^N \|\mathbb{S}(:,j) - \sum_{k=1}^n (\mathbb{S}(:,j),\mathbb{M}(:,k))_{\mathbb{R}^{N_h}} \mathbb{M}(:,k) \|_{\mathbb{R}^{N_h}}^2.
\end{equation*}
In conclusion, a \emph{POD-Galerkin} ROM consists of solving \eqref{eq:system_ROM} using the matrix $\mathbb{V}=\widetilde{\mathbb{W}}(\nrb)$ to define the reduced space.

%%%%%%%%%%%%%%%%%%%%%%%%%%%%%%%%%%%%%%%%%%%%%%%%%%%%%%%%%%%%%%%%%%%%%%
\subsection{Limitations of linear methods and nonlinear model order reduction}\label{sec:benchmark}

Despite mathematically rigorous, POD-Galerkin ROMs might feature some limitations when applied to problems with microstructure, as shown by some simple numerical test cases that we shall discuss in this section. As a prototype problem, we consider \eqref{eq:problem} with an operator independent on the parameters, \ie $L_{\para} = -\Delta$, with homogeneous Dirichlet boundary conditions, and the action of the microstructure entirely carried out by the right hand side $f_{\para}=f_{\parmi}$. We consider a unit square domain $\Omega=(-1,1)^2$ partitioned using a 50$\times$50 computational mesh $\mathcal{T}_h$ of uniform triangular elements. For the sake of simplicity, we discretize the problem using piecewise linear, continuous, Lagrangian finite elements $V_h=X_h^1(\Omega)$, and  consider two representative instances of the following characteristics of the microstructure: {\em (i)}  a problem with continuously variable scales, and {\em (ii)}  a problem with scale separation.
 
\paragraph{Problem with continuously variable scales} To model this feature of the microstructure, we consider 
\[
f_{\parmi}(x,y) = p(x) p(y), \ \ \mbox{with} \ \  p(z) =  \sum_{k=1}^n (\alpha_k^z \sin(k \pi z) + \beta_k^z \cos(k \pi z))
\]
and $\parmi = \{\alpha_k^z, \beta_k^z\}_{k=1}^n$ with $n=6$ and $z=x,y$, for a total of 24 parameters modulating the linear combination of 144 trigonometric functions. For the application of the POD-Galerkin method to this test case, we sample the parameter space assuming a uniform probability distribution $\mathcal{U}[-1,1]$ for all parameter, selecting a total amount of $N=1000$ snapshots to build the data matrix $\mathbb{S}$. Equivalent results are obtained for a larger data matrix $\mathbb{S}$ made of 2000 snapshots-- therefore, we keep $N=1000$ as a reference value for the forthcoming tests.

\begin{figure}
    \includegraphics[height=0.2\textheight]{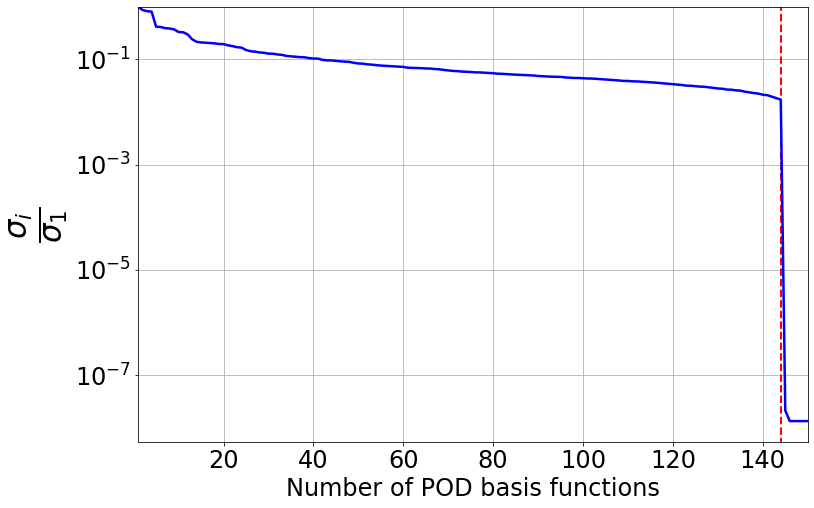}
    \includegraphics[height=0.2\textheight]{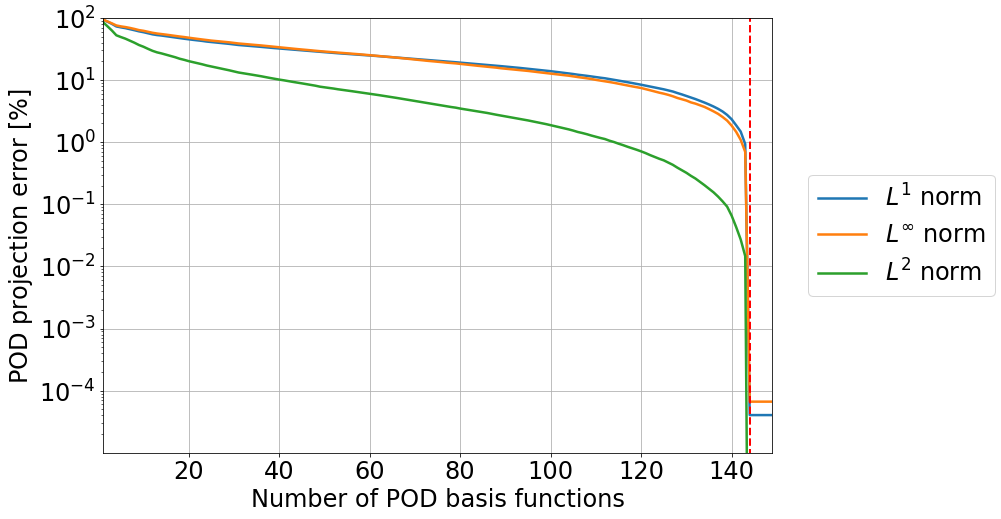}
    \caption{Decay of the singular values of $\mathbb{S}$ and the convergence rate of the projection error with respect to $\nrb$ for a problem with continuously variable scales. The red line in the left panel denotes the threshold of $n=144$ POD modes.}
    \label{fig:toy1}
\end{figure}

Since we are interested to assess the capability of the POD basis to approximate the solution manifold, in this case we do not calculate the actual solution of \eqref{eq:system_ROM}, rather we consider the projection error obtained when projecting the FOM solution onto the POD basis for an increasing dimension $\nrb$ of the reduced basis, that is, 
\begin{equation}\label{eq:EPOD}
    E_{POD}(\nrb;\uuph) = \frac{\| (\mathbb{I}-\VV(\nrb)\VV^T(\nrb))\uuph \|_{2,N_h}}{\| \uuph \|_{2,N_h}}
\end{equation}
where $\VV(\nrb)$ denotes the matrix encoding the linear approximation space spanned by $\nrb$ POD basis functions and   $\| \mathbf{v} \|_{p,N} = \big(\sum_{i=1}^{N}|v_i|^p\big)^{1/p}$ is the $p$-norm in the generic vector space $\mathbb{R}^{N}$. 
We note that the solution of this problem depends linearly on the parameters, as a result if we apply the POD-Galerkin method with at least 144 basis functions, we shall represent exactly the discrete solution manifold -- in other words, $E_{POD}(\nrb=144;\uuph)=0$ for any $\uuph \in \RR^{N_h}$. Actually, the desired behavior of the ROM would be to achieve a satisfactory approximation with much less than 144 basis functions. 
The decay of the singular values of $\mathbb{S}$ and the convergence rate of the projection error with respect to $\nrb$, measured for a parameter value and corresponding finite element function not included in the snapshot matrix, are reported in Figure \ref{fig:toy1}. 

These results immediately show that, due to the slow decay of the singular values, the performance of the POD-Galerkin method cannot be satisfactory in this case.

\paragraph{Problem with scale separation} In this second test case, we discuss the behavior of the POD-Galerkin method for a problem where the forcing term is made of the superimposition of trigonometric functions with different frequencies, with a gap in the frequency spectrum. In the same computational setting described before, we consider the function 
\[
f(x,y) = f_{\parma}(x,y)+f_{\parmi}(x,y), \ \ \mbox{with} \ \  f_{\parma}(x,y)=p_L(x)p_L(y),\, f_{\parmi}(x,y)=p_H(x)p_H(y)
\]
and
\begin{equation*}
    P_L(z) = \sum_{k=1}^2 (\alpha_k^z \sin(k \pi z) + \beta_k^z \cos(k \pi z)),
    \quad
    P_H(z) = \sum_{k=5}^8 (\widetilde{\alpha}_k^z \sin(k \pi z) + \widetilde{\beta}_k^z \cos(k \pi z)),
\end{equation*}
for a total of 8 coefficients encoding the low frequency scales and 16 parameters for the high frequency scales, all sampled using a uniform probability function $\alpha_k^z, \beta_k^z, \widetilde{\alpha}_k^z, \widetilde{\beta}_k^z \sim \mathcal{U}[-1,1])$. 
These parameters modulate the forcing term made by the linear combination of 80 trigonometric modes.

\begin{figure}
    \includegraphics[height=0.2\textheight]{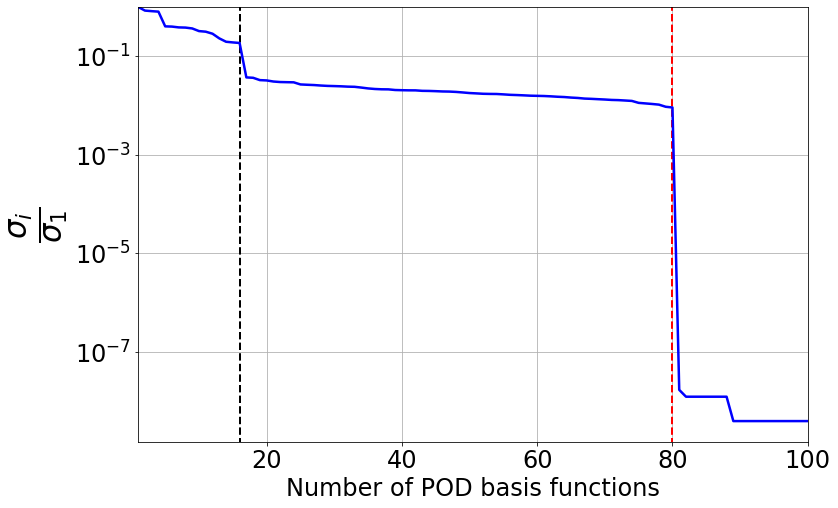}
    \includegraphics[height=0.2\textheight]{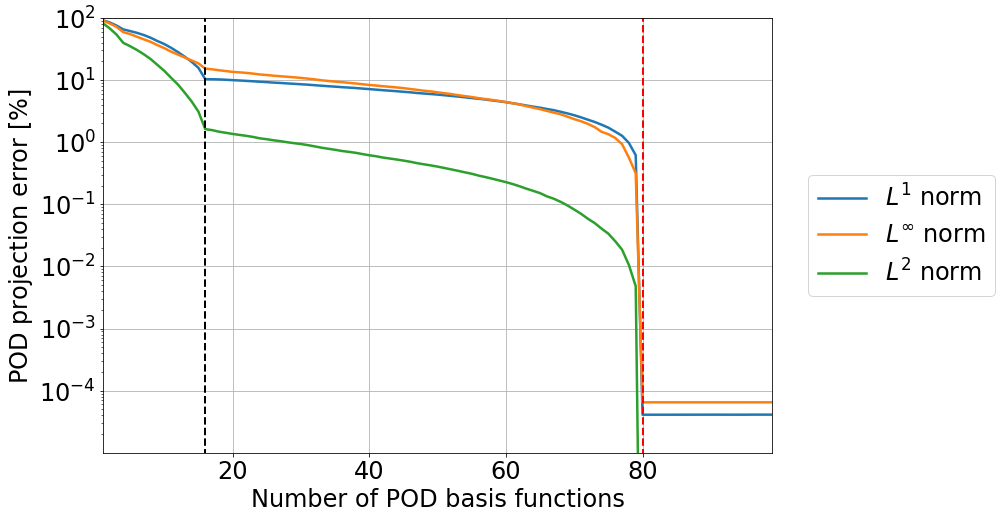}
    \caption{Decay of the singular values of $\mathbb{S}$ and the convergence rate of the projection error with respect to $\nrb$ for a problem with scale separation. The vertical lines in the left panel denote the number of parameters of the low frequency function ($n=16)$ and the total number of parameters ($n=80$).}
    \label{fig:toy2}
\end{figure}

The purpose of this test case is to investigate whether the POD basis with a number of entries comparable to the dimension of the low frequency source terms can approximate well the whole solution space. Figure \ref{fig:toy2} shows that the decay of the singular values reflects the gap in the frequency spectrum. However, after a reasonably good performance of the POD method in approximating the effect of the low frequencies on the solution, a large region of stagnation of the approximation properties can be highlighted.

Therefore, it is mandatory to include the effect of the microscale in the ROM approximation space, with a subsequent increase of the number of POD modes required to achieve a satisfactory approximation of the FOM.

% \paragraph{Problem with unknown scales}
% In this case, the forcing term is affected by two single frequencies, one of which is unknown. Precisely, for the same problem and the same discretization settings described in the previous tests, we consider $f(x,y) = \sin(\delta \pi x) (\alpha \sin(\pi y) + \beta \cos(\pi y))$, with uniformly distributed coefficients $\alpha, \beta \sim U([-1,1])$ and variable frequency $\delta$, randomly chosen in the interval $\delta \sim U([2, 10])$. The function $f$ has only three degrees of freedom, but the dependence of $f$ on $\delta$ is \emph{nonlinear}. For this reason, there is no guarantee that the projection error will be small just a few (approximately three) modes. 

% This conjecture is confirmed by the results of Figure \ref{fig:toy3}, where we report the decay of the singular values of $\mathbb{S}$ and the decay of the projection error with \nrb, the number of basis functions. We see that the approximation properties of the POD basis become barely satisfactory only when a number of basis $n>>3$ is considered.

% \begin{figure}
%     \includegraphics[height=0.2\textheight]{Images/singularvaluesDOD.png}
%     \includegraphics[height=0.2\textheight]{Images/errprojDOD.png}
%     \caption{Decay of the singular values of $\mathbb{S}$ and the convergence rate of the projection error with respect to \nrb for a problem with unknown scales.}
%     \label{fig:toy3}
% \end{figure}

\subsubsection{Non intrusive, nonlinear reduced order modeling}
The previous examples show that a Kolmogorov barrier arises due to the presence of the microstructure, this latter limiting the decay of the error achieved with projection-based ROMs  that seek linear approximations in spaces. The main objective of this work is to overcome these limitations resorting to nonlinear model order reduction, implemented by combining a linear ROM approximation with nonlinear maps built through deep feedforward NNs, thus introducing a nonlinear, and non intrusive, ROM. 

Nonlinear ROMs encompass a wide class of strategies aiming at replacing the linear approximation hypothesis encoded in \eqref{eq:lin_reduced_basis} with more general, nonlinear maps. Rather than replacing the linear combination of POD modes with fully nonlinear approximations, as done, e.g., through deep learning-based ROMs \cite{FrancoROMPDE,fresca2021comprehensive,fresca2021pod}, we express our approximation following the approach proposed in \cite{Peherstorfer2022725}, as 
\begin{equation}\label{eq:nonlin_reduced_basis}
    \nuphrb = \sum_{i=1}^{\nnrb} \nurb^{(i)}\npsi_i(x,\nft(\para)),
\end{equation}
where the main difference with respect to \eqref{eq:lin_reduced_basis} is that the basis functions $\npsi_i$ now depend on the parameters of the problem through some \emph{features}, named $\nft(\para)$. 
Stated differently, the nonlinear approximation \eqref{eq:nonlin_reduced_basis} is a linear combination of functions $\npsi_i(x,\nft(\para))$, $i=1,\ldots, \nnrb$ that depend through $\nft(\para)$ on the element $\uph$ being approximated; this is indeed different from the linear approximation \eqref{eq:lin_reduced_basis}, where the basis functions are fixed and independent of which element $\uph$ of the (discrete) solution manifold is approximated. As we will see in the following, proceeding in this way can be extremely helpful in order to break the Kolmogorov barrier, leading to lower errors  for the same number of degrees of freedom.

For the particular case of problems with microstructure, we put into action this general nonlinear ROM framework by using the nonlinear features as an additive correction, named \emph{closure}, to the linear approximation obtained through \eqref{eq:lin_reduced_basis}. This consists in seeking a reduced approximation of the FOM based on the following representation,
    \begin{equation}\label{eq:corrected_reduced_basis}
        \nuphrb = \sum_{i=1}^{\nrb} \urb^{(i)}\psi_i(x) + \nurb\npsi(x,\nft(\para)),
    \end{equation}
where the function $\npsi(x,\nft(\para))$ will be learned through artificial NNs (introduced in detail in the next section), mapping the parameter space into the full order FOM space, namely $\ann: \mathcal{P} \to \RR^{N_h} \equiv V_h$. In other terms, the linear approximation space is replaced by $\widetilde{V}_{rb} = V_{rb} + \ann(\mathcal{P})$. 

The presence of the microstructure poses however an additional challenge to the one of the Kolmogorov barrier discussed so far. Indeed, space-varying parameters deserve a similar dimensionality reduction as done for the state solution, given their dimension comparable to the number of high-fidelity FOM degrees of freedom. Moreover, the potentially involved nature of space-varying parameters makes the use of classical hyper reduction techniques -- that would employ, e.g., a linear POD basis to approximate those data in order to comply with the ROM requirements -- not feasible. For this reason, we resort to the use of NN-based approximations, exploiting a particular class of networks, the so-called \emph{mesh informed neural networks} (MINNs) previously proposed by the authors \cite{FrancoMINN}, to handle mesh based functional data in a very efficient way; the resulting strategy, involving a POD-based dimensionality reduction for the state solution, and a MINN-based approximation of space-varying fields, will be referred to as POD-MINN method. A mesh informed architecture is also employed to approximate the ROM coordinates, similarly to the POD-NN approach \cite{hesthaven2018non}, but extended here to the case of MINNs to handle the high-dimensional input data. The final approximation, obtained when the closure model is applied to the POD-MINN approximation, will be referred to as POD-MINN+. %, see Figure \ref{fig:POD-MINN+}.

In the next two sections we describe the formulation, the implementation details and the performance of these methodologies applied to problems with microstructure. Then, we will apply them to a real-life problem modeling that, despite some simplifying assumptions, is capable of representing the effect on microcirculation on the tissue microenvironment, involving multiple spatial scales and reaching very tiny detail levels.

%%%%%%%%%%%%%%%%%%%%%%%%%%%%%%%%%%%%%%%%%%%%%%%%%%%%%%%%%%%%%%%%%%%%%%
\section{Mesh-informed neural networks} %Artificial neural networks in the approximation of PDEs}
\label{sec:nn} 

\textit{Artificial neural networks} (ANN) are parametrized functions between two (typically high-dimensional) vector spaces that have been recently employed also for the approximation of PDEs in several contexts. In this section, we discuss how they can be designed and exploited in particular case of problems with microstructure, to obtain non-intrusive and accurate ROMs. Before addressing mesh-informed neural networks, representing a key tool of our methodology, we briefly review feedforward neural networks for the sake of notation.

\subsection{The architecture and the training of feedforward neural networks}

ANNs are computational models obtained by stacking compositions of nonlinear transformations, acting on a collection of nodes called \textit{neurons}, that are organized into building blocks called \textit{layers}, these latter linked together through weighted connections. Given two vector spaces $V_m \subset \mathbb{R}^m$ and $\RR^n$,  a layer $\lnn$ that takes a vector in $\RR^m$ as input and uses the \textit{activation function} $\rho$, is a map of the form
\begin{equation*}
    \lnn: \RR^m \to \RR^n \ \textrm{such that} \
    \lnn(\mathbf{v}) = \rho(\WW\mathbf{v} + \mathbf{b}), 
\end{equation*}
where $\WW\in \mathbb{R}^{n\times m}$ and $\mathbf{b} \in \RR^m$ are the \emph{weight matrix} and the bias of the layer, respectively. The activation function $\rho$, to be selected, is applied to the linear combination componentwise for each neuron of the same layer and ensures the nonlinearity of the approximation. 
The topology of the connections between the neurons determines the architecture of the ANN. The most common and known example is the \emph{feedforward neural network}, defined as the composition of multiple layers $\lnn_i\,:\RR^{n_{i}}\rightarrow \RR^{n_{i+1}}$, with $i=1,\ldots,l$, where $\lnn_{l+1}$ is the output layer and the remaining ones are called \textit{hidden layers}. A feedforward neural network produces a map of the form $\ann:=\lnn_{l+1}\circ \lnn_{l}\ldots \circ \lnn_1$, so that, denoting with $\WW_i$ the weight matrices between the i-th and the $i+1$-th layer, $\mathbf{b}_i$ the corresponding bias and $\rho$ the activation function (generally the same for each layer):
\begin{equation*}
    \ann(\mathbf{v}; \WW_1, \dots, \WW_{l+1}, \mathbf{b}_1, \dots, \mathbf{b}_{l+1}) =\rho_{l+1}(\WW_{l+1}\rho_{l}(\ldots\WW_2  \rho_1(\WW_1\mathbf{v}+\mathbf{b}_1)+\mathbf{b}_2\ldots)+\mathbf{b}_{l+1}).
\end{equation*}
For the sake of clarity, we make a distinction between the notation related to an NN seen as an operator -- that is, a parametrized function denoted with $\ann$ and depending on the input vector $\mathbf{v} \in \RR^{n_1}$, the weights and biases $\WW_i,\mathbf{b}_i$ for $i=1,\ldots,l+1$ --  from the output of the same neural network, that is, a vector $\mathbf{w}\in\RR^{n_{l+1}}$. To lighten the notation, we denote the hyperparameters of the neural network, namely the collection of weights and biases, with the general symbol $\hyper$. 
The choice of the activation function is left to the user and is problem-dependent. Furthermore, it can be neglected in the output layer \cite{Kutnyok, Schwab}.

% \begin{figure}[p]
%     \centering
%     \includegraphics[width=0.8\textwidth]{Images/NN_plain.png}
%     \caption{Schematic representation of a feedforward neural network. Notice that the biases are not reported in order to prevent cluttering.}
%     \label{fig:NN_plain}
% \end{figure}

After determining the architecture of the neural network, the hyperparameters for the training process are set and tuned. Hinging upon the choice of the activation function for each layer, the initial state is provided: weights and biases are set to zero or randomly initialized in order to retain the variance of the activation across every layer, avoiding the exponential decrease or increase of the magnitude of input signals. 
In the context of \textit{supervised learning}, the training of a neural network consists of tuning its weights and biases, simply put $\hyper$, by means of the minimization of a suitable \textit{loss} function $\mathcal{E}$ that measures the discrepancy between a given dataset and the predictions of the neural network. The minimization is typically performed by means of a gradient descent algorithm \cite{Schmidhuber}, such as the L-BFGS method \cite{LiuNocedal} (a quasi-Newton optimizer) or ADAM \cite{Adam}.

The computational costs of the training are heavily linked to the \textit{complexity} of the neural network, that depends on its depth (the number of layers), on the number of neurons within each layer and on possible constraints on the weighted and biases. Regarding the latter, the architecture is \textit{dense} if no constraints are imposed. In the following subsection, we discuss a different approach, that introduces sparsity in the weights matrix to efficiently approximate FE functions, exploiting their mesh-based information. 

\subsection{Mesh-informed neural networks for the approximation of finite element functions}

%\todo[inline]{L'applicazione di POD-NN a problemi con microstruttura richiede di usare dati in input ad elevata dimensionalità definiti sulla mesh, quindi è necessario usare le MINNs}

Despite their incredible expressivity, plain NN architectures based on dense layers can incur into major limitations, especially when dealing with high-dimensional data. In fact, when the dimensions into play become fairly large, dense architectures quickly become intractable, as they result in models that are harder to optimize and often prone to overfitting (thus also requiring more data for their training). Unfortunately, within our context, these issues arise quite naturally. In fact, the geometrical features of the microstructure are encoded at a high-fidelity level: consequently, any architecture working with the microstructure information must be powerful enough to handle an input of dimension $O(N_{h})$, with $N_{h}\gg 1$.

To overcome this difficulty, we rely on a particular class of sparse architectures called Mesh-Informed Neural Networks (MINNs), originally proposed in \cite{FrancoMINN} as a way to handle high-dimensional data generated by FE  discretizations. MINNs are obtained through an a priori pruning strategy that promotes local operations rather than global ones. We can summarize the general idea underlying MINNs as follows. Let $\mathbb{R}^{N_{h}}\cong V_{h}$ and $\mathbb{R}^{M_{h'}}\cong V_{h'}$ be two vector spaces, each associated to a given FE discretization. The two may refer to different discretizations in terms of mesh step size or polynomial degree, but they must be defined over a common spatial domain $\Omega$. 
\newcommand{\x}{\mathbf{x}} 

Let $$\{\x_{j}\}_{j=1}^{N_{h}}\subset\Omega\quad\text{and}\quad\{\x_{i}'\}_{i=1}^{M_{h'}}\subset\Omega,$$ be the coordinates of the nodes for the two FE spaces, respectively. Then, following the definition in \cite{FrancoMINN}, a mesh-informed layer with support $r>0$ is a DNN layer $\lnn:\RR^{N_{h}}\to\RR^{M_{h'}}$ whose weight matrix $\WW$ satisfies the sparsity constraint below
\begin{equation*}
	\WW_{i,j}=0\quad\forall i,j,\;\text{such that}\;\;d(\x_{j},\x_{i}')>r,
\end{equation*}
where $d:\Omega\times\Omega\to[0,+\infty)$ is a suitable distance function defined over the spatial domain, such as, e.g., the Euclidean distance  $d(\x,\x'):=|\x-\x'|$ -- this choice, however, is not restrictive, making then possible to use also a geodesic distance. Here, the intuition is that each entry component $j$ at input (resp. output) has been associated to some node $\x_{j}$ representing a degree of freedom in the FEspace at input (resp. output). Then, this interpretation allows one to construct a sparse DNN model whose weight matrix ultimately acts as a local operator in terms of the corresponding FE spaces. In general, for simplicial meshes in $\Omega\subset\mathbb{R}^{d}$, it can be shown that the weight matrix of a mesh-informed layer with support $r>0$ has at most $O(r^{d}N_{h}M_{h'})$ nonzero entries \cite{FrancoMINN}. The advantage of such construction is that it reduces the number of hyparameters to be optimized during the training stage without affecting the overall expressivity of the model. 

In particular, training a MINN model can be far less demanding than doing the same with its dense counterpart; furthermore, the sparsity constraints naturally introduced in a MINN can have a beneficial impact on the generalization capabilities of the architecture, see, e.g., \cite{FrancoMINN} for further details. 

In light of these considerations, for what concerns our purposes, we shall exploit both mesh-informed and dense layers: the choice will depend, from time to time, on the nature of the data at hand.
To distinguish between MINNs and general dense neural networks, we denote by $\minn$ those architectures that feature the use of mesh-informed layers. 
%\textcolor{blue}{Vero: vogliamo usare $\mathcal{S}$ per il solution manifold? Btw, ho rifrasato: denoterei con $\mathcal{M}$ architetture anche solo "parzialmente mesh-informed", visto che $\mathcal{M}_{rb}$ è ibrida (parte micro MINN e macro densa)}\textcolor{green}{PZ: correzione per solution manifold fatta}
%For further details about MINNs and their implementation, we refer the reader to \cite{FrancoMINN}.

%\subsection{Some results for the appoximation of PDEs using neural networks}

%%%%%%%%%%%%%%%%%%%%%%%%%%%%%%%%%%%%%%%%%%%%%%%%%%%%%%%%%%%%%%%%%%%%%%
\section{Non-intrusive and nonlinear ROMs: POD-MINN and POD-MINN+ strategies}\label{sec:pod-nn}

In this section we illustrate a strategy to use mesh-informed neural networks for the construction of non-intrusive ROMs for problems with microstructure, detailing the approximation introduced in equation \eqref{eq:corrected_reduced_basis}. We remark that we will exploit MINNs %neural networks 
in two distinct steps:
\begin{enumerate}
\item in the first step, directly inspired to \cite{hesthaven2018non}, MINNs are used to approximate the map from the parameters of the problem to the reduced coefficients using the network 
\[\minn_{rb}: 
\quad \mathcal{P} \rightarrow \RR^{\nrb},
\quad \para \mapsto \muurb; \]
%\textcolor{red}{Però questo è in contraddizione con il diagramma, in cui compare anche $\mathcal{N}_{rb}$. Oltretutto la MINN dovrebbe solo "ridurre" $\mu_m$, mentre $\mu_M$ dovrebbe entrare in gioco allo stesso livello dell'approssimazione  MINN di $\mu_m$.}\textcolor{green}{Modifico il diagramma per renderlo piu coerente con la formula (14)}
\item in the second step, MINNs are exploited for the \emph{closure} model, namely a nonlinear correction to the previous ROM (based on a linear approximation space). To this purpose, we introduce the MINN $\minn_c$ such that
%\[ \ann_c: 
%\quad \mathcal{P} \rightarrow \RR^{N_h},
%\quad \para \mapsto \nurb\widetilde{\psi}(x,\widetilde{\alpha}(\para)).\]
\[ \minn_c: 
\quad \mathcal{P} \rightarrow \RR^{N_h},
\quad \para \mapsto \nuuu.\]
\end{enumerate}
%\textcolor{red}{AM: Si rischia di fare confusione tra cosa è linear e cosa è nonlinear. Io intendo con linear ROM un modello ridotto di tipo POD-Galerkin. Un modello POD-MINN (o anche solo POD-NN) è linear? Secondo me no. Abbiamo bisogno di definirlo tale? Perché non distinguere tra intrusive e non-intrusive, piuttosto?
%}\\\textcolor{blue}{NF: My two cents: POD-NN e POD-MINN sono "nonlinear" nella stessa misura in cui POD-Galerkin lo è. In quest'ultimo infatti la soluzione è della forma $u\approx \mathbf{V}\text{solve}(\boldsymbol{\mu})$, dove $\text{solve}(\boldsymbol{\mu})$ è la risoluzione del problema proiettato (per ogni $\boldsymbol{\mu}$ fissato). In generale, la mappa $\boldsymbol{\mu}\mapsto\text{solve}(\boldsymbol{\mu})$ sarà sempre nonlineare, anche per POD-Galerkin a coefficienti affini. In questo senso, sia POD-Galerkin che POD-MINN sono dei: linear ROM se interpretiamo la linearità come espansione sui modi; nonlinear ROM se parliamo di dipendenza dai parametri (ma a quel punto qualunque ROM lo sarebbe). Se ci sembra confusivo, possiamo evitare il tipo di nomenclatura (non usarei intrusive vs nonintrusive, visto che sia POD-MINN che la chiusura sono nonintrusivi).}
%\\\textcolor{green}{PZ: sono d'accordo con NF}

\begin{figure}[b!]
    \centering
    \includegraphics[width=\textwidth]{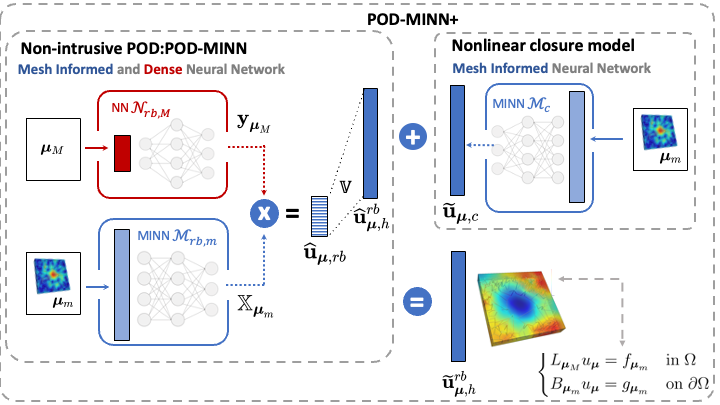}
    \caption{A sketch of the POD-MINN+ method. The macroscale parameters, $\parma$, and the microscale ones, $\parmi$, are fed to two separate architectures,  $\ann_{rb,M}$ and $\minn_{rb,m}$, whose outputs are later combined to approximate the RB coefficients, cf. Eq. \eqref{eq:POD_NN_split}. The coefficients are then expanded over the POD basis, $\mathbb{V}$, and the ROM solution is further corrected with a closure term computed by a third network, $\minn_{c}$, that accounts for the small scales.} 
    \label{fig:POD-MINN+}
\end{figure}

The overall framework is summarized in  Figure \ref{fig:POD-MINN+}. We recall that parameters affecting the differential problem feature a {\em scale separation property}, with macroscale parameters showing a small dimensional and, conversely, microscale parameters that are instead high-dimensional -- that is, $\npmicro\gg \npmacro$. Approximating quantities depending on the former can thus be done through plain feedforward neural networks, while we decide to encode the latter through a mesh-informed neural network to approximate spatially-dependent parametric fields related to the microscale. We shall detail how in the upcoming subsection.

%%%%%%%%%%%%%%%%%%%%%%%%%%%%%%%%%%%%%%%%%%%%%%%%%%%%%%%%%%%%%%%%%%%%%%
\subsection{Approximation of the %reduced basis
POD coefficients using mesh informed neural networks} 

We observe that the best approximation of the FOM solution $\uph$ -- whose degrees of freedom are collected in the vector $\uuph \in \RR^{N_h}$ -- in the POD space is given by its projection 
\[
\VV \VV^T \uuph = \sum_{i=1}^{\nrb} \big(\VV^T \uuph\big)^{(i)}\psi_i(x).
\]
The expression above defines a mapping 
\begin{equation}\label{eq:pi}
\pi: 
\quad \mathcal{P} \subset \RR^p \rightarrow \RR^{\nrb},
\quad \para \mapsto \VV^T \uuph,
\end{equation}
whose approximation could provide a way to generate, onto the POD space, the solution of the problem for unseen values of the input parameters. As shown in \cite{hesthaven2018non}, if the dimension of the reduced space is small enough, it is possible to effectively approximate the map $\pi$ by means of a neural network. Since for the case of our interest the parameter space $\mathcal{P}$ includes the high-dimensional data necessary to encode the microscale parameters $\parmi$, we use a MINN instead of a generic feedforward neural network.

%\textcolor{red}{Qui occorrerebbe secondo me essere più chiari sulla distinzione tra la MINN per la microscala e la NN per la macroscala, se le usiamo entrambe.}

We remark that this method overrides the projection step of the FOM onto the reduced space, \eqref{eq:projection_FOM}, and the solution of the reduced problem \eqref{eq:system_ROM}. As a consequence, it turns out to be a completely non-intrusive reduced order modeling strategy that will be particularly advantageous to handle problems that do not enjoy an affine decomposition of the FOM operators, or in which one is interested to approximate only a subset of the variables involved in its formulation.

The dataset to train the network $\minn_{rb}$ is the representation in the reduced space of the FOM solutions collected into the data matrix $\SS$. In other words, the training set are the input-output couples $\{\para^{(i)},\VV^T\SS(:,i)\}_{i=1}^N$, where $\para^{(i)}$ are the parameter values used to calculate the FOM solution vector $\SS(:,i)=\uuphi=\uuh(\para^{(i)})$. This dataset is used to define the loss function,
\begin{equation*}
    \mathcal{E}(\minn_{rb};\hyper_{rb}) = \frac{1}{N} \sum_{i=1}^{N} \| \VV^T\SS(:,i) - \minn_{rb}(\para^{(i)};\hyper_{rb})\|_{2,\nrb},
\end{equation*}
and the (perfectly) trained neural network is the one characterized by the paramaters $\hyper_{rb}^*$ that minimize the loss, precisely $\hyper_{rb}^*=\argmin_{\hyper} \mathcal{E}(\minn_{rb};\hyper_{rb})$. 
% In general the training algorithms do not converge toward the global minimum. The error due to convergence toward local minima of the loss function is called optimization error.

%\textcolor{red}{AM: Detto così sembra che gli stessi dati usati per costurire lo spazio POD siano usati una seconda volta... non si prendono nuovi snapshots?} \textcolor{blue}{NF: credo sia siano usati gli stessi (in fase POD, si guardano solo le soluzioni e si "butta via" l'informazione sui parametri, che qui viene ripresa in considerazione).}\textcolor{green}{PZ: esatto}

Once the network has been trained, the model can be used in an online phase similarly to the baseline reduced basis approach: given a new set of parameters $\para$, the mapping through the network $\muurb = \minn_{rb}(\para;\hyper_{rb}^*)$ and a subsequent linear transformation, $\VV \muurb$, enable the approximation of the FOM solution $\uuph(\para)$.
The approximation provided by $\ann_{rb}$ can be represented by the following function,
\begin{equation}\label{eq:POD_NN}
    \muphrb(x) 
    = \sum_{i=1}^{\nrb} \muurb^{(i)}\psi_i(x)
    = \sum_{j=1}^{N_h} {\muuphrbj}\phi_j(x),
\end{equation}
where $\muuphrb = \VV \muurb \in \RR^{N_h}$ is the vector collecting its degrees of freedom in the FEM space.
\\\\
The POD-MINN method is based on specific architectures designed with the purpose to manage differently the parameters related to the \emph{micro} and \emph{macro} scales, as they belong to spaces with different dimension. The first step, $\minn_{rb}$ combines a MINN that manages the high-dimensional information about $\parmi$ with a DNN that deals with $\parma$, leading to the following formulation of an input-output map that approximates the reduced coefficients:
\begin{equation}\label{eq:POD_NN_split}
\minn_{rb}(\para; \hyper_{rb}) = \minn_{rb,\micro}(\parmi; \hypermi)\ann_{rb,\macro}(\parma; \hyperma),
\end{equation}
where:
\begin{itemize}
    \item The output of $\ann_{rb,\macro}: \RR^{\npmacro} \rightarrow \RR^{k}$ is a vector $\mathbf{y}_{\parma}=\ann_{rb,\macro}(\parma) \in \RR^{k}$;
    \item The network $\minn_{rb,\micro} : \RR^{\npmicro} \rightarrow \RR^{n_{rb} \times k}$ gives back a matrix $\mathbb{X}_{\parmi} \in \RR^{n_{rb} \times k}$. 
\end{itemize}
We note that the network $\minn_{rb}$ approximates the mapping $\pi$ defined in \eqref{eq:pi}. 
%In this particular architecture, we use hyperbolic tangent activation functions and the approximation of the reduced coefficients provided by $\minn_{rb}$ is represented by a linear combination of vectors $\mathbb{X}(:,j)$ and hyper-coefficients $\mathbf{y}_j,\,j=1,\ldots,k$.
%In some sense, the matrix $\mathbb{X}_{\parmi}$ represents local bases for the reduced coefficients, $\muurb$, that depends on the microscale. 
The final output is then represented by 
\[\muurb = \mathbb{X}_{\parmi}\cdot\mathbf{y}_{\parma}, %= \sum_{j=1}^k \mathbb{X}(:,j)\mathbf{y}_j.
\]
where $\cdot$ emphasizes the presence of a matrix-vector multiplication. 
The training of the network defined in \eqref{eq:POD_NN_split} is done all at once using the loss function $\mathcal{E}(\minn_{rb};\hyper_{rb})$.
The respective weights and biases, collected in the vectors $\hypermi$ and $\hyperma$, are drawn from the normalized initialization proposed in \cite{GlorotBengio}.

Although the approach proposed here provides an efficient way to recover the main features of the FOM solution, it is characterised by two main error sources, in fact, regardless of the training of the network, by orthogonality of $\VV$:
\begin{equation*}
    \|\uuph - \VV\muurb\|_{2,N_h}^{2}=\|\uuph - \VV\VV^{T}\uuph\|_{2,N_h}^{2}+\|\VV^{T}\uuph-\muurb\|_{2,\nrb}^{2}.\end{equation*}
This is to emphasize that the neural network only provides an estimate of the map $\pi$, meaning that $\muurb \approx \VV^{T}\uuph$; furthermore, on top of this, the representation of the FOM solution in the reduced space is affected by a projection error $\|\uuph - \VV\VV^{T}\uuph\|_{2,N_h}$, identifying a lower bound for the approximation of the FOM solution that can not be overcome by the POD-MINN approach. The approximation of problems with microstructure entails the choice of a higher number of reduced basis functions in order to keep the projection error sufficiently small and, consequently, the reconstruction error of the solution. We tackle these issues by means of a closure model, as shown in the following section.

%%%%%%%%%%%%%%%%%%%%%%%%%%%%%%%%%%%%%%%%%%%%%%%%%%%%%%%%%%%%%%%%%%%%%%
\subsection{Mesh informed neural networks for the closure model}

Designing a closure model for problems with microstructure is particularly challenging for the following reasons:
\emph{(i)} the input of the network $\minn_c$ must necessarily include the microscale parameters $\parmi$ that are represented in the high-dimensional FE space relative to the FOM model;
\emph{(ii)} the output of the closure model is also isomorphic to the FE space of the FOM;
\emph{(iii)} as a consequence of \emph{(ii)}, the training data of $\minn_c$ are also high-dimensional.
For these reasons, MINNs play a crucial role for the successful application of the closure model.

Let us denote with $\nuuu$ the output of the closure model, namely $\nuuu=\minn_c(\para) \in \RR^{N_h}$. Then, following \eqref{eq:corrected_reduced_basis}, the closure model is easily added to the reduced basis approximation,
\begin{equation}\label{eq:POD_NN_closure}
    \nuphrb(x) = \sum_{j=1}^{N_h} \nuuu^{(j)}\phi_j(x) 
    + \sum_{i=1}^{\nrb} \muurb^{(i)}\psi_i(x),
\end{equation}
where $\psi_i(x)$ are the reduced bases while $\phi_j(x)$ are the standard FE basis functions used in the FOM model.
In vector form, the previous equation is equivalent to 
\[
\nuuphrb = \nuuu + \VV\muurb.
\]
The dataset for training the closure model consists of the FOM snapshots collected in the data matrix $\SS$, more precisely the couples $\{\para^{(i)},\SS(:,i)\}_{i=1}^N$. 
Moreover, we make the choice to feed the closure model only with the microscale input parameters. Precisely, we have
\begin{equation}
    \minn_{c}(\para; \hyper_c) = \minn_c (\parmi; \hyper_c).
\end{equation}
As a result, we look for a network $\minn_c (\parmi; \hyper_c): \RR^{\npmicro} \to \RR^{N_h}$ that minimizes the following loss function,
\begin{equation*}
    \mathcal{E}(\minn_c;\hyper_c) = \frac{1}{N} \sum_{i=1}^{N} \| \big(\SS(:,i) - \VV\minn_{rb}(\para^{(i)};\hyper_{rb}^*)\big)
    - \minn_c(\parmi^{(i)};\hyper_c) \|_{2,N_h},
\end{equation*}
where $\SS(:,i) - \VV\minn_{rb}(\para^{(i)};\hyper_{rb}^*)$ represents the approximation error generated by the POD-MINN method, related to the $i^{th}$-snapshot, which is then corrected by $\minn_c$ to further reduce it. The weights and biases of $\hyper_c$ are initialized to zero, since the optimization process starts from an initial state determined by the POD-MINN reconstruction of the solution.

%The output of $\minn_c$ is denoted by $\mathbf{z}_{\parmi}=\minn_c (\parmi) \in \RR^{N_h}$. 
As result, the POD-MINN+ approximation of the FOM solution is given by the following expression:
\begin{equation*}
  \widetilde{\mathbf{u}}_{\para,h}^{rb} 
  =\VV \muurb + \nuuu
  =\VV\mathbb{X}_{\parmi}\mathbf{y}_{\parma} + \minn_c (\parmi) \approx \uuph.
\end{equation*}

In conclusion, augmenting the POD-MINN trained map through the closure model allows to retrieve also the variability of local features in the solution manifold that POD-MINN can recover only considering a large number of reduced basis functions, including the small frequency POD modes. The method with closure is called POD-MINN+.

\subsection{Numerical results obtained with the  POD-MINN and POD-MINN+ methods}

In this section we analyze the performance of the POD-MINN+ method in the approximation of the benchmark problems proposed in section \ref{sec:benchmark}. 
To assess the performance of the closure model, we compare the POD-MINN+ reconstruction error with the POD-MINN one and the projection error generated projecting the FOM solution onto the POD space, defined as in \eqref{eq:EPOD}. The former two are defined as follows:
\begin{gather}\label{eq:EPODMINN}
    E_{\tiny POD-MINN}(\nrb,\uuph)=\frac{\|\uuph - \VV(\nrb)\muurb \|_{2,N_h}}{\| \uuph \|_{2,N_h}},
\\
    E_{\tiny POD-MINN+}(\nrb,\uuph)=\frac{\|\uuph - \VV(\nrb)\muurb - \nuuu\|_{2,N_h}}{\| \uuph \|_{2,N_h}}.
\end{gather}

For both benchmark problems addressed in Sect.~\ref{sec:benchmark}, we consider the same architectures for the POD-MINN method and the closure model, hinging upon the use of MINNs. The microscale input parameters $\parmi$ are the representations of the respective forcing terms $f(x,y)$ in the finite-element space of the solution $V_h$, with dimension $N_h$. Let $\mathcal{T}_h$ denote the computational mesh (collection of all geometric elements) related to the space $V_h$.
We remind that we consider a unit square domain $\Omega=(-1,1)^2$ partitioned using uniform triangles corresponding to a 50$\times$50 grid named  $\mathcal{T}_h$.

As a first step, related to the POD-MINN approach, we introduce the map approximating the reduced coefficients defining the architecture $\minn_{rb} (\parmi; \hyper_{rb}): \RR^{N_h} \to \RR^{n_{rb}}$, exploiting a mesh-informed layer $L^{rb}_1$ of support $r = 0.6$ and a dense layer $L^{rb}_2$:
$$L^{rb}_1: \,\, V_h \xrightarrow{r=0.6} V_{2h},\quad\quad L^{rb}_2: V_{2h} \xrightarrow{\textcolor{white}{r=0.6}} \RR^{n_{rb}}$$
%\begin{itemize}
%    \item $L^{rb}_1$: $\,\, V_h \xrightarrow{r=0.6} V_{2h}$, 
%    \item $L^{rb}_2$: $V_{2h} \xrightarrow{\textcolor{white}{r=0.6}} \RR^{n_{rb}}$,
%\end{itemize}
with $V_{2h}=X_{2h}^1(\Omega)$ finite-element space defined over the coarser computational mesh $\mathcal{T}_{2h}$ (stepsize $2h$, 25$\times$25 grid), made by uniform triangles, of the same domain $\Omega$. 
For the closure model, let $\mathcal{T}_{h^\prime}$ a coarse mesh in $\Omega$, with ${h^\prime} > h$. We use a MINN architecture $\minn_c (\parmi; \hyper_c): \RR^{N_h} \to \RR^{N_h}$, composed by 2 mesh-informed layers  $L^c_i$ of support $r=0.6$:
$$L^c_1 : \, V_h \xrightarrow{r=0.6} V_{h^\prime},\quad\quad L^c_2 : V_{h^\prime} \xrightarrow{r=0.6} V_h$$
%\begin{itemize}
%    \item $L^c_1 : \, V_h \xrightarrow{r=0.6} %V_{h^\prime}$,
    %\item $L^c_2 : V_{h^\prime} \xrightarrow{r=0.6} V_h$,
%\end{itemize}
where $V_{h^\prime}=X_{h^\prime}^1(\Omega)$ is the finite-element space defined from a 35$\times$35 computational mesh $\mathcal{T}_{h^\prime}$ of uniform triangles. 
As result, the MINNs $\minn_{rb} (\parmi; \hyper_{rb})$ and $\minn_c (\parmi; \hyper_c)$ are defined as follows:
%\begin{flalign*}
$$    %&
\minn_{rb} : V_h \xrightarrow{r=0.6} V_{2h} \xrightarrow{\textcolor{white}{r=0.6}} \RR^{n_{rb}}, \quad\quad%\\
    %&
    \minn_c : \,\, V_h \xrightarrow{r=0.6} V_{h^\prime} \, \xrightarrow{r=0.6} V_{h}.
$$%\end{flalign*}
For what concerns the training of the neural networks, the $N=1000$ snapshot functions, taken into account to build the data matrix $\mathbb{S}$, are partitioned to provide:
\begin{itemize}
    \item a \textit{training} set of $N_{train}=750$ samples;
    \item a \textit{validation} set consisting of $N_{valid}=50$ snapshots;
    \item a \textit{test} set of $N_{valid}=200$ snapshots.
\end{itemize}
The optimization of the loss function is performed through the L-BFGS optimizer with learning rate equal to 1, without batching. The networks are trained for a total of $250$ epochs, using an early stopping criterion based on the validation error, that is applied if the following conditions are met: the training error \textit{decreases} but the validation error \textit{increases} for at least two consecutive epochs.
%\begin{itemize}
    %\item the training error \textit{decreases} over the last 2 epochs,
    %\item the validation error \textit{increases} in the last 2 epochs.
%\end{itemize}

The training times for the POD-MINN+ method in the two benchmark problems vary between from 98 to 130 seconds for the test case with continuous scales and from 89 to 188 seconds for the test case with continuously variables scales. As a general trend, we denote that an increase of the number of POD basis functions entails a decrease of the training time of the closure model.

% \begin{table}[htb]
% \begin{center}
%  \begin{tabular}{c|ccc}
%  \textbf{Benchmark problem} & \textbf{Minimum} $[s]$ & \textbf{Maximum} $[s]$ & \textbf{Average} $[s]$ \\
%  \hline
%  Continuosly variable scales & $89$ & $188$ & $158$ \\
%  Scale separation  & $98$ & $130$ & $117$ \\
% \end{tabular}   
% \end{center}
% \caption{Training times for POD-MINN+ for the two case studies
% reported in Section 2.5. For each benchmark problem, we apply POD-MINN+ varying the number of POD basis functions and we report the minimum, the maximum and the average training time.}
% \label{tab::trainingtimesbenchmark}
% \end{table}

% \begin{figure}
%     \centering
%     \includegraphics[height=0.17\textheight]{Images/PODprojectionerrorsPODclosureFULL.png}
%     \includegraphics[height=0.17\textheight]{Images/PODprojectionerrorsPODclosureLOWHIGH.png}
%     \caption{The performance of the POD-MINN+ method, namely $E_{\tiny PODMINN+}$ defined in \eqref{eq:EPODMINN+}, compared to the one of the POD projection defined as $E_{POD}$ in \eqref{eq:EPOD}, measured in the norms $p=1,2,\infty$. The case with continuously variable scales is reported on the left, the one with the scale separation is on the right.}
%     \label{fig:PODMINN+}
% \end{figure}

\begin{figure}
\centering
\includegraphics[height=0.19\textheight]{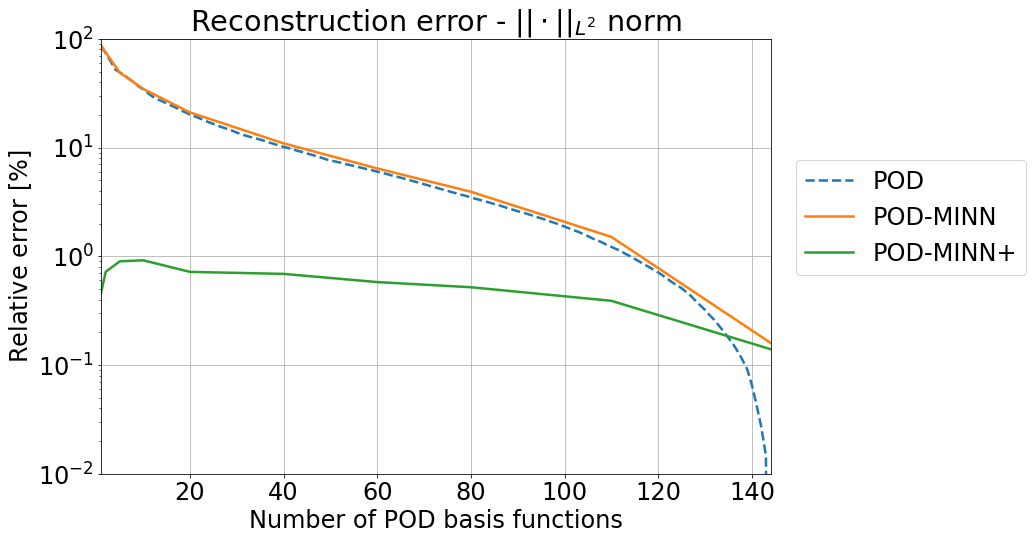}
\includegraphics[height=0.19\textheight]{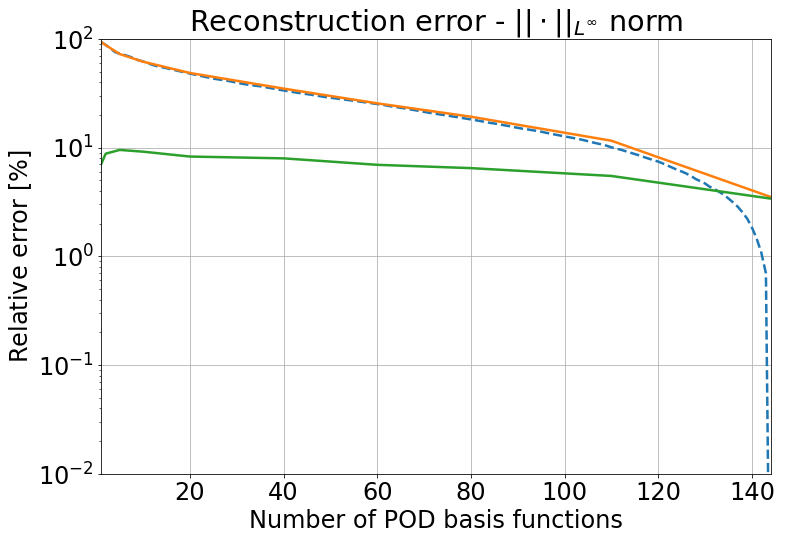}
\\
\includegraphics[height=0.19\textheight]{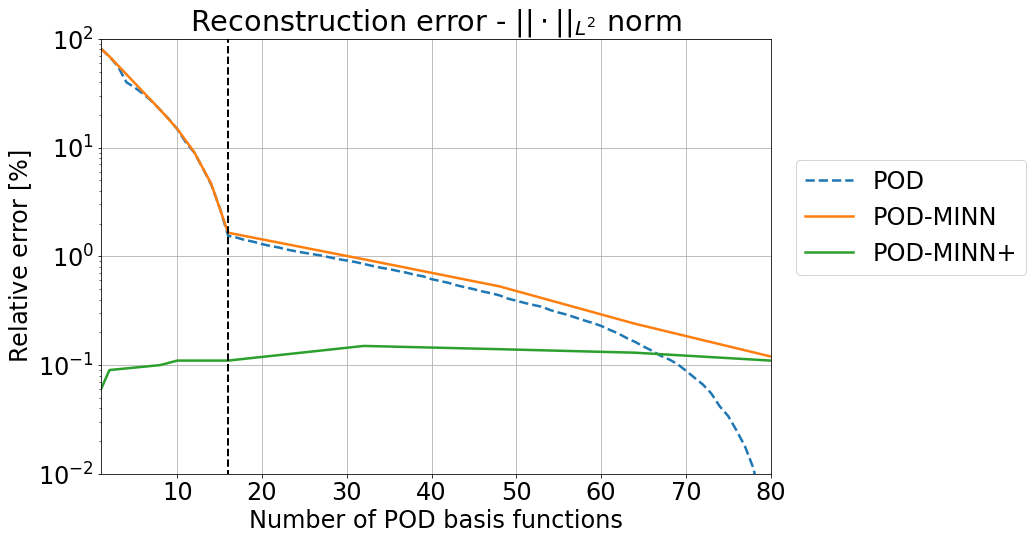}
\includegraphics[height=0.19\textheight]{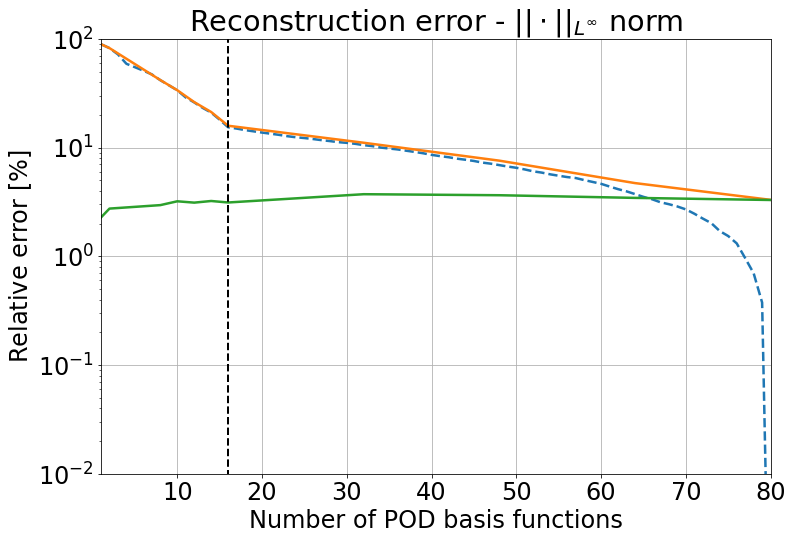}
\caption{The performance of the POD-MINN+ method applied to the benchmark problems with continuous scales (top row) and separate scales (bottom row). We compare the errors $E_{\tiny POD-MINN}$ and $E_{\tiny POD-MINN+}$ defined in \eqref{eq:EPODMINN}, to the one of the POD projection defined as $E_{POD}$ in \eqref{eq:EPOD}, measured in the norms $p=2,\infty$. The case with $p=2$ is reported on the left, the one with $p=\infty$ is on the right.}
\label{fig:benchmarkPODMINN+}
\end{figure}

%\todo[inline]{Aggiornare questa discussione con il confronto tra POD, POD-MINN e POD-MINN+ come definiti sopra.}

The comparison between the POD projection error, namely $E_{POD}$, and $E_{\tiny PODMINN+}$ when varying the number of basis functions $\nrb$ is reported in Figure \ref{fig:benchmarkPODMINN+} for the benchmark problems with continuously variable scales and with scale separation. These results confirm that the the POD-MINN approach is not performing well, because it can not provide a better approximation than the POD projection error. The POD-MINN error is indeed very close to the POD projection, but this approximation is qualitatively not satisfactory for both the $L^2$ and $L^\infty$ norms, unless a large number of bases is used.

Conversely, the POD-MINN+ method is very effective especially in the regime with a low number of basis functions. In particular, for the case with continuously variable scales the POD-MINN+ method reduces of 2 orders of magnitude the POD error in the region $\nrb < 20$, when it is measured in the 2-norm. For the benchmark problem with the scale separation the gain increases to more than three orders of magnitude, see Figure \ref{fig:benchmarkPODMINN+} (bottom-left panel). This effect is particularly evident when we approximate this benchmark problem with basis functions only able to capture the low frequency modes, namely the case $\nrb \leq 16$. In this regime the closure model is entirely responsible for the approximation of the high frequency modes. The POD-MINN+ method effectively performs this task with a relative error of $0.1 \% = 10^{-3}$.

% Overall we see that in the regime of number of basis functions significantly smaller than the number of modes characterizing the solution manifold, the closure model breaks the Kolmogorov barrier in the model reduction phase.

We note that the closure model does not improve the rate of convergence with respect to $\nrb$, it only decreases the magnitude of the error. In general, the error of the POD-Galerkin method enjoys an optimal exponential decay with respect to $\nrb$ that can not be increased resorting to the closure model. This observation confirms that the POD-MINN+ method is meaningful for those problems where the Kolmogorov $n$-width and consequently the POD projection error decays slowly with respect to $\nrb$.

\section{Application to oxygen transfer in the microcirculation}
\label{sec:apps}

In this section we present an application of the POD-MINN+ method to model oxygen transfer at the level of microcirculation, described by \eqref{eq:oxy}. After introducing the FOM model and its discretization, we discuss its parametrization and finally we present the numerical results that illustrate the advantages of the proposed reduced order model.

\subsection{Numerical approximation of oxygen transfer in the microcirculation: the full order model}

For the spatial approximation of the problem we use the finite-element method. We consider a 3D domain $\Omega$, identified by a slab of edge $1\,mm$ and thickness $0.15\,mm$ and discretized through a $20\times20\times3$ structured computational mesh of tetrahedra. For the interstitial domain we introduce the space of the piecewise linear, continuous, Lagrangian finite-elements $V_{t,h}=X_h^1(\Omega)$, with dimension $N_h = 1764$. 
On the other hand, we assume that the 1D vascular network is immersed in the 3D slab and we discretize it by partitioning each vascular branch $\Lambda_i$ into a sufficiently large number of linear segments. A piecewise linear, continuous, Lagrangian finite-element space $V_{v,h}^i=X_h^1(\Lambda_i)$ is employed for each branch $\Lambda_i$, where the index $i$ spans through the vascular branches $i=1,\ldots,N_b$. Hence, we approximate the equations in the whole microvascular network through a finite element space $V_{v,h}=\Big(\bigcup_{i=1}^{N_b}V_{v,h}^i \Big) \bigcap C^0(\Lambda)$.

\begin{figure}[htb]
\centering
\includegraphics[height=0.22\textheight]{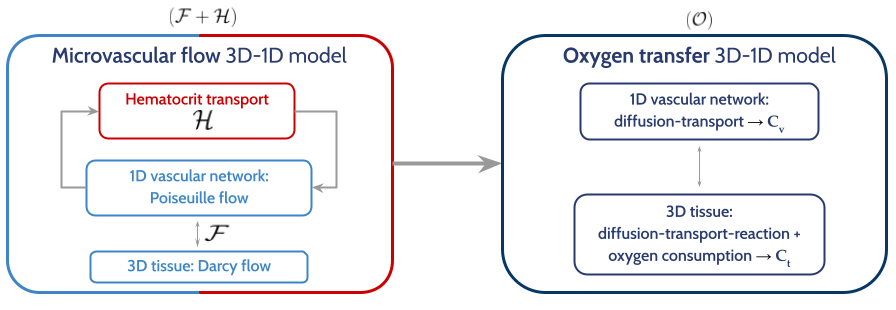}
\caption{General layout of the full order model for the whole vascular microenvironment.}
\label{fig:layoutFOM}
\end{figure}

It is important to clarify that the finite-element approximation of the oxygen transfer model is performed downstream to the full comprehensive problem for the vascular microenvironment, described in \cite{Possenti2019,Possenti20213356}. As it is shown in Figure \ref{fig:layoutFOM}, indeed, we firstly compute via finite-element method the velocity and the pressure in the tissue and the vascular network, together with the discharge hematocrit in the vascular network. We remark that this proposed ROM based on the POD-MINN+ method could have been equivalently applied to any variable in the domain $\Omega$, such as the pressure or the velocity fields.

% \subsection{Output of the problem}

\subsection{The parametrization of the FOM: inputs and outputs}

Since the model presented in Section \ref{sec:setup} is driven by parametrized PDEs, a crucial aspect of the problem is the choice of which parameters to include as inputs of the  model. We consider the parameter space $\mathcal{P} = \mathcal{P}_{\macro}\cup\mathcal{P}_{\micro}$, where the macro- and micro-scale subspaces are selected according to the results of an earlier sensitivity analysis \cite{VitulloSensitivity}. We then select  three physical parameters (see Table \ref{tab:parameters})   among those appearing in the model \eqref{eq:oxy}, and two geometrical, mesh-based inputs, encoding the information related to the 1D vascular networks. The range of variation of each input parameter is determined consistently with respect to its physiological bounds.

\begin{figure}[htb]
    \centering
    \includegraphics[scale=0.16]{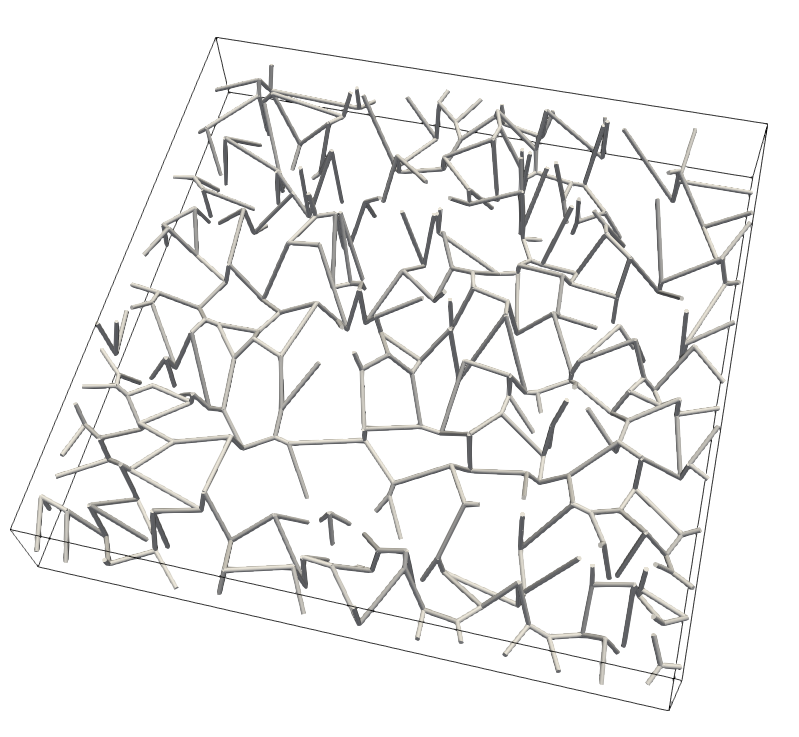}
    \includegraphics[scale=0.16]{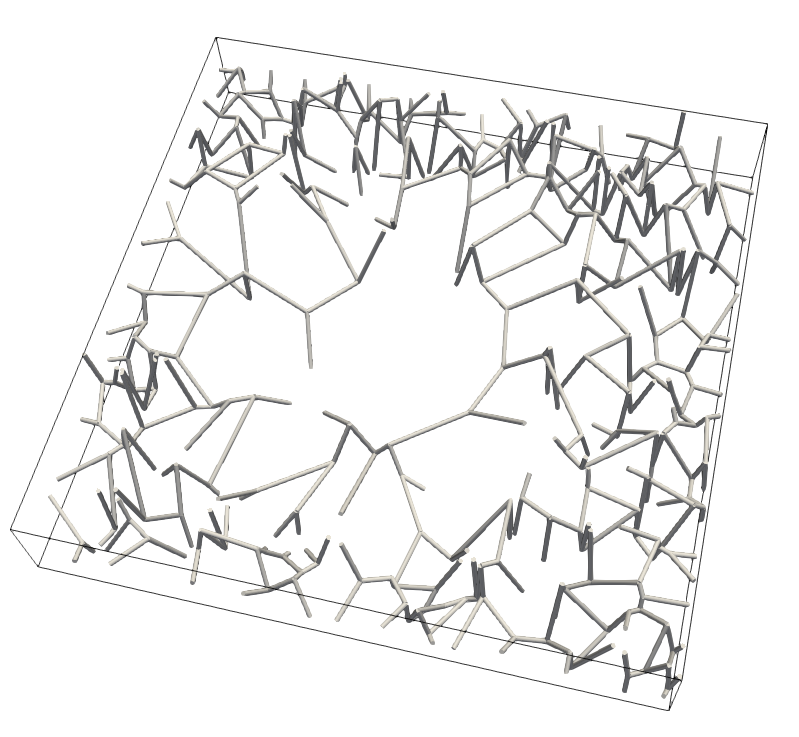}
    \includegraphics[scale=0.16]{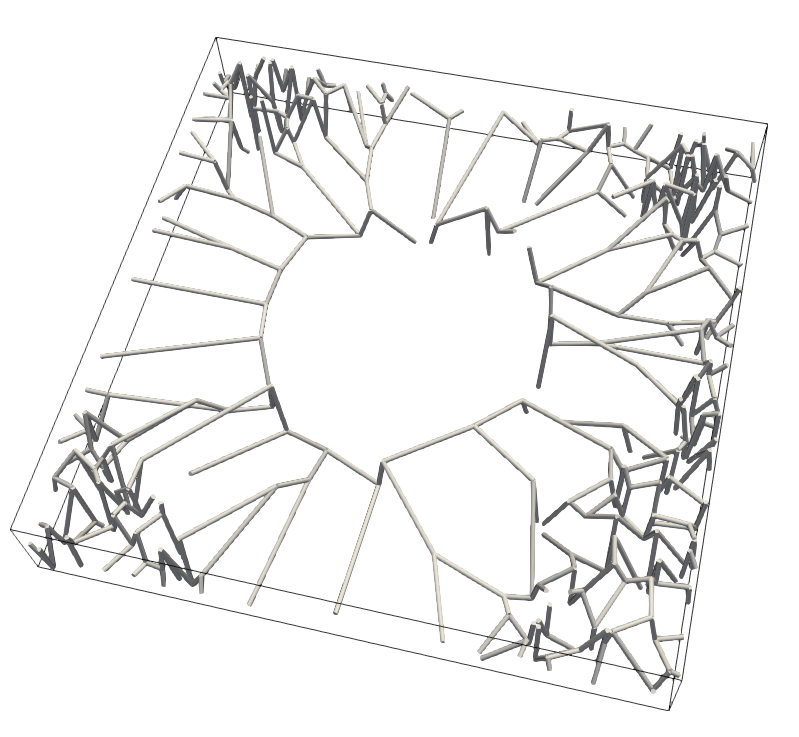}
    \caption{Examples of artificial vascular networks with extravascular distance $\textbf{d}$ progressively increasing from left to right.}
    \label{fig:networks}
\end{figure}

In particular, the vascular architectures are obtained starting from a biomimetic algorithm that replicate the essential traits of angiogenesis, starting from some parameters that characterize the density and the distribution of small vessels in a vascularized tissue. These quantities are the vascular surface per unit volume of tissue, named $S/V$, and an indicator that governs the distribution of point seeds used to initialize the angiogenesis algorithm. The definitions and the range of these hyper-parameters of the vascular model are reported in Table \ref{tab:parameters}.
By sampling these hyper-parameters and applying the angiogenesis algorithm, we obtain a population of admissible vascular networks with sufficient variability. Figure \ref{fig:networks} shows some examples of 1D vascular networks with different spatial distributions, spanning from regular structures to less ordered ones. For more details regarding the generation of the embedded microstructures, we refer to \cite{VitulloSensitivity}. From the computational standpoint, the vascular networks are represented as metric graphs, the vertices of which are points in the 3D space. This description is not practical for the purpose of generating a reduced order model. It is more convenient to transform the vascular graph into two continuous functions. The main one is the distance function from the nearest point of the vascular network, defined on $\Omega_t$ and named $\mathbf{d}$. The second one is used to identify the intersections of the vascular network with the boundary and it is named $\boldsymbol \eta$. These are the input functions of the ROM and are presented below.

\subsubsection{Extravascular distance and inlet function}
% The full order model exploits a numerical decoupling to allow for a separate treatment of the solution of the intravascular flow and that of the interstitial flow: due to this, an independent choice for the grids and discretizations of the capillary net and the tissue slab can be made, since the relationship between the two solutions is encoded in the mathematical formulation of the problem. 
As a way to retrieve an equivalent information of the 1D vascular graph we thus choose to introduce the \textit{extravascular distance} \ie the function mapping each point in the tissue slab to its distance from the closest point of the vascular network. We represent this function as a linear and piecewise continuous finite element function in $V_{t,h}$. The vector of degrees of freedom of such discrete function, corresponding to its nodal values, is  $\dist \in \RR^{N_h}$.
The vector $\dist$ plays the role of the forcing term $f_{\parmi}$ in the model \eqref{eq:micro_problem}.

\begin{figure}[htb]
\centering
\includegraphics[height=0.21\textheight]{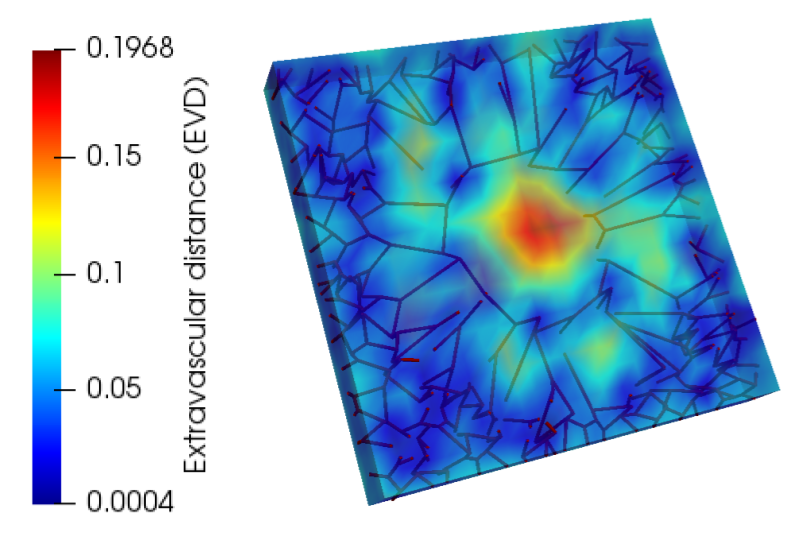}
\includegraphics[height=0.21\textheight]{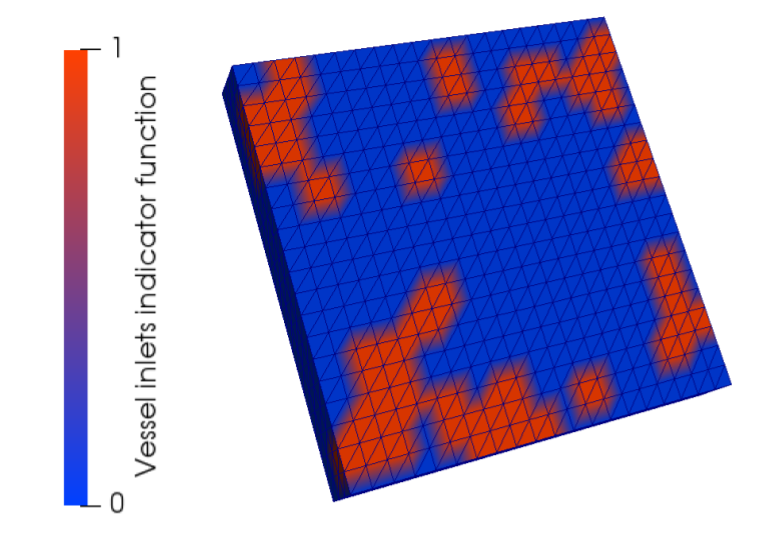}
\caption{The nodal distance from the nearest vessel of the vascular network is shown on the left. On the right we plot the indicator function of the vessel inlets on the boundary.}
\label{fig:inputmicro}
\end{figure}

Furthermore, it is useful to provide the ROM with the information of the inlet points of the vascular network on the boundary. To this purpose, we introduce an \textit{inlet characteristic function} that assigns a unit value to nodes of the tissue slab grid which are close to the inlets. Again, this function is defined at the discrete level using the space $V_{t,h}$. In practice, a vertex of the finite element mesh is marked as having non-zero value if and only if one of the surface elements is intersected by an inlet vessel. The degrees of freedom of this function are denoted by $\inlets \in \RR^{N_h}$. Figure \ref{fig:inputmicro} shows the mesh-based data $\dist$ and $\inlets$ for a particular instance of the vascular graph.

We point out, without loss of generality, that we can extend this approach encoding the geometry of the vascular network with high-dimensional data defined over coarser meshes (\textit{super-resolution}) or refined ones (\textit{sub-resolution}) \cite{Fukami2023}. 

% \subsubsection{Inlets function - $\inlets$}
% While the entire 1D graph can be seen as a source term for the diffusion problem in the interstitial tissue, it is clear that these will be stronger in proximity of capillaries on which inlet conditions have been taken compared to those with outflow conditions. Indeed, the oxygen concentration must have decreased as a consequence to the permeability of the vascular walls. To recover this asymmetry in the vascular bed, we introduce an \textit{inlets characteristic function} $\inlets \in \RR^{N_h}$, which assigns a unit value to nodes of the tissue slab grid which are close to the inlets: in particular, a vertex is marked as having non-zero value if and only if one of the surface elements it is part of is intersected by an inlet capillary.

\subsubsection{Output of the full order model}
\begin{figure}[htb]
\centering
\includegraphics[height=0.19\textheight]{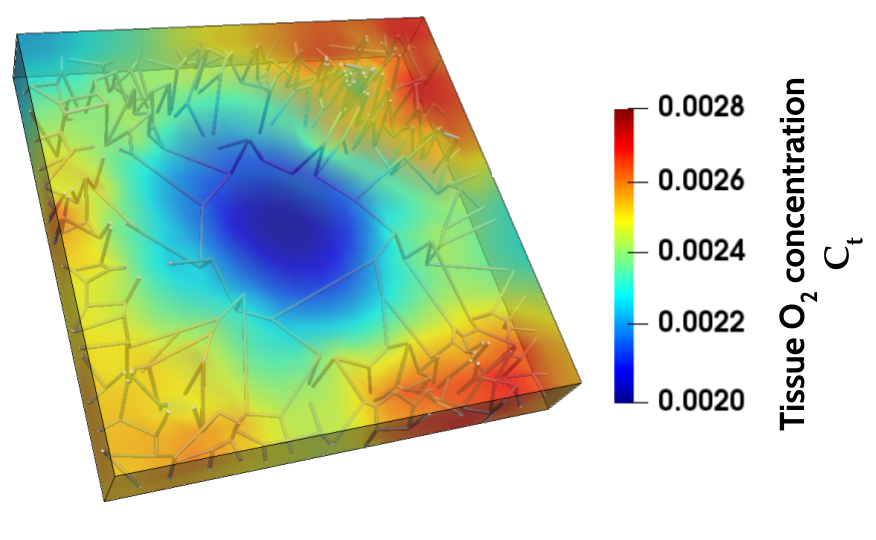}
\includegraphics[height=0.19\textheight]{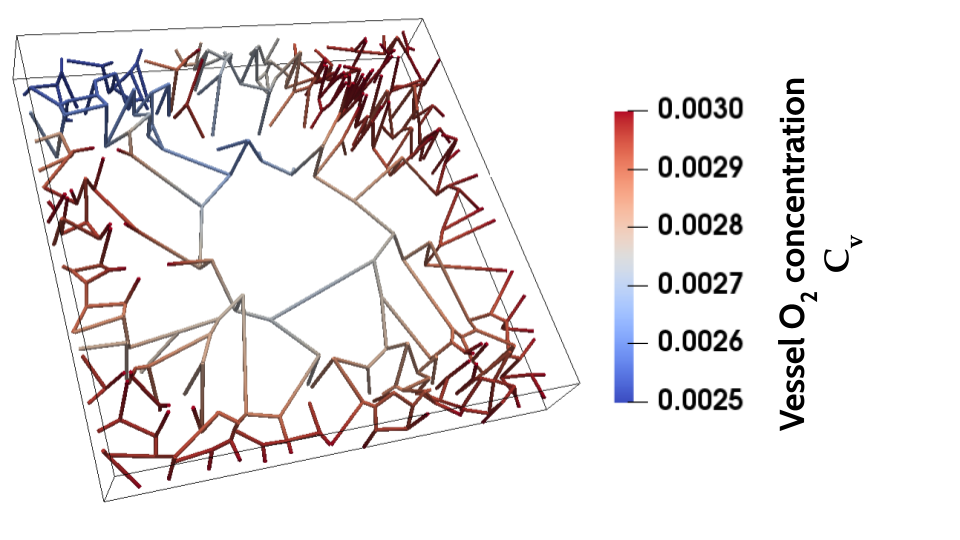}
\caption{The FOM solution representing the tissue oxygenation map ($mL_{O_2}/mL_B$) is reported on the left panel. On the right we show the 1D embedded vascular microstructure that visibly influences the oxygen map, over which the vessel oxygen concentration $C_v$ is represented.}
\label{fig:outputmicro}
\end{figure}

The output of the FOM is the tissue oxygenation map $C_t$, one of the state variables of the oxygen transfer model, measured in $mL_{O_2}/mL_B$. In the equations \eqref{eq:micro_problem} it represents the high-fidelity solution $u_{\para}$. An example of the oxygen map $C_t$ is shown in Figure \ref{fig:outputmicro}, where we notice that the embedded vascular microstructure has a significant influence of on the tissue oxygenation. We also note that for the development of the ROM the variable $C_t$ is suitably rescaled by a factor $300$. We remark that we can exploit the 3D-1D full order model of the vascular microenvironment to extend the proposed approach for the approximation of other state variables such as the interstitial pressure in the tissue domain $\Omega_t$. On the contrary, the methodology has to be revised if the aim is to reconstruct outputs defined over the vascular network.

 \begin{table} [h!tb]
	\begin{center}
		\begin{tabular}{||c c c c||}
		\hline
			\textsc{symbol} & \textsc{Parameter} & \textsc{Unit} & \textsc{Range of variation} \\ \hline
			$P_{O_2}$ & \small{$O_2$ wall permeability} & $m/s$ & $3.5 \times 10^{-5}-3.0 \times 10^{-4}$
            \\[6pt]
			$V_{max}$ & \small{$O_2$ consumption rate} & $mL_{O_2}/cm^3/s$ & $4.0 \times 10^{-5}-2.4 \times 10^{-4}$
            \\[6pt]
		$C_{v,in}$ & \small{$O_2$ concentration at the inlets} & $mL_{O_2}/mL_B$ & $2.25 \times 10^{-3}-3.75 \times 10^{-3}$
            \\[6pt]
   		$\% \frac{SEEDS_{(-)}}{SEEDS_{(+)}}$ & \small{Seeds for angiogenesis} & [$\%$] & $0-75$ 
        \\ [1ex] 
            $S/V$ & \small{Vascular surface per unit volume} & [$m^{-1}$] & $5\cdot10^3-7\cdot10^3$ \\ [1ex] \hline
		\end{tabular}	
            \caption{The first three rows illustrate the biophysical parameters of the ROM and their ranges of variation. The last two rows report the hyper-parameters used to initialize the algorithm that generates the vascular network.}
      	\label{tab:parameters}
	\end{center}
\end{table}

\subsection{Implementation of the POD-MINN and POD-MINN+ methods}

\begin{figure}[htb]
\centering
\includegraphics[height=0.21\textheight]{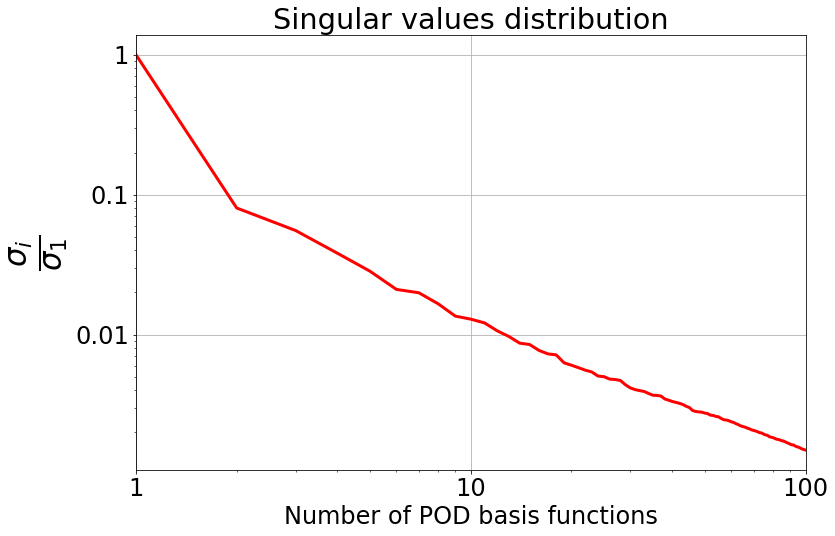}
\includegraphics[height=0.21\textheight]{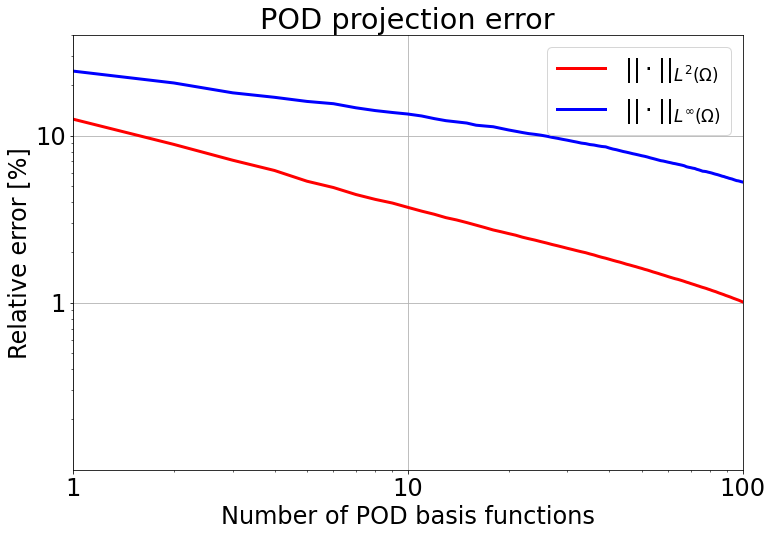} 
\caption{(\textit{Left}) The decay of the singular values distribution for the microcirculation problem (\textit{Right}) The POD projection error of the FOM solution on the reduced basis space varying the number of basis functions.}
\label{fig:poderrors}
\end{figure}

Once defined the macroscale and the microscale parameters, namely $\parma$ $= [V_{max},C_{v,in},P_{O_2}]$ and $\parmi$  $= (\dist, \inlets)$, we apply a Monte Carlo sampling of the parameter space and  collect the corresponding input-output pairs provided by the FOM model, named $\{(\parma^{(i)}, \dist^{(i)}, \inlets^{(i)});  \uuph^{(i)} \}_{i=1}^N$. We organize these data into the data matrix  $\SS$. In particular, for what concerns the microscale parameter space, the sampling of the vascular networks is performed with respect to the hyper-parameters introduced in Table \ref{tab:parameters}. Thanks to the chosen settings, the extravascular distance $\dist$ and the inlets function $\inlets$ share the same dimension and ordering of the snapshots in the matrix $\SS$, so that the matching between inputs and outputs can be easily established.
The entries of the physical parameters $\parma$ are normalized between 0 and 1 with respect to their range of variation reported in Table \ref{tab:parameters}.
The available dataset of $N = 1600$ FOM snapshots and parameters samples is partitioned as follows:
\begin{itemize}
\item $N_{train}=1220$ training data; 
\item $N_{valid}=100$ validation data;
\item $N_{test}=280$ testing data.
\end{itemize}
We perform a singular value decomposition of the training dataset, building the projection matrix $\VV$, collecting the discrete representation of $\nrb$ basis functions. As a preliminary analysis, we study in Figure \ref{fig:poderrors} the convergence rate of the projection error of the FOM output on the reduced space. These results immediately show that, since the problem is globally diffusion-driven and overall well described by the macroscale features at the low frequencies, the decay of the singular values is relatively quick, meaning that a POD approach is expected to yield good results, as shown in Figure \ref{fig:poderrors} (\textit{left}). Nonetheless, by comparing projection errors measured in the the $L^2$ and the $L^{\infty}$ norms, reported in Figure \ref{fig:poderrors} (\textit{right}), we see that the $L^\infty$ error dominates over the $L^2$ one, meaning that the local effects due to the microstructure still have a big impact on the whole solution. Rather than including a large number of POD modes, it is therefore preferable to retain just few of them, and then include a suitable closure model,
%useful to rely on a  \textcolor{red}{nonlinear} ROM, 
with the aim to retrieve the local features of the oxygen map and improve the overall accuracy of the method.

\subsubsection{Architectures and training of the POD-MINN method}

In order to implement the POD-MINN method, we exploit two independent MINNs to process the mesh-based geometrical input data $\dist$ and $\inlets$, using the same architectures in both cases.

We recall the formulation of the input-output function $\minn_{rb}$ in (\ref{eq:POD_NN_split}):
\begin{equation}\label{eq:POD_MINN_micro}
\minn_{rb}(\para; \hyper_{rb}) = \minn_{rb,\micro}(\parmi; \hypermi)\ann_{rb,\macro}(\parma; \hyperma),
\end{equation}
where, in this case, the MINN that takes as input the geometrical parametrization is defined as
\begin{equation*}
\minn_{rb,\micro}(\parmi; \hypermi) = \minn_{rb,\eta}(\inlets; \hyperi) \odot  \minn_{rb,d}(\dist; \hyperd), 
\end{equation*}
splitting the contributions of the single geometrical input parameters, with the inlets function $\inlets$ weighting the effect of the extravascular distance $\dist$ through the MINN $\minn_{rb,\eta}$. Here, we denote with $\odot$ the Hadamard product between two matrices,
\begin{equation*}
    (A \odot B)_{i, j} = (A)_{i,j} (B)_{i,j} \qquad \forall A, B \in \mathbb{R}^{m \times n}, \; m,n \geq 1.
\end{equation*}

The neural network $\minn_{rb,d}(\dist; \hyperd): \RR^{N_h} \to \RR^{n_{rb}\times k}$ is assembled relying on two mesh-informed layers $L^{rb,d}_1$ and $L^{rb,d}_2$ of support $r = 0.3$ and hyperbolic tangent activation function, combined with a dense layer $L^{rb,d}_3$ complemented with a suitable reshape of the output, which results in the following structure:
%\begin{itemize}
%    \item $L^{rb,d}_1$: $\,\, V_{t,h} \xrightarrow{r=0.3} V_{t,H}$,
%    \item $L^{rb,d}_2$: $\,\, V_{t,H} \xrightarrow{r=0.3} V_{t,H}$,
%    \item $L^{rb,d}_3$: $V_{t,H} \xrightarrow{\textcolor{white}{r=0.3}} \RR^{n_{rb}k}$,
%    \item Reshape: $\RR^{n_{rb}k} \xrightarrow{\textcolor{white}{r=0.3}} \RR^{n_{rb}\times k}$,
%\end{itemize}
%$$L^{rb,d}_1: \,\, V_{t,h} \xrightarrow{r=0.3} V_{t,H},\quad\quad L^{rb,d}_2: \,\, V_{t,H} \xrightarrow{r=0.3} V_{t,H}$$
%$$L^{rb,d}_3: V_{t,H} \xrightarrow{\textcolor{white}{r=0.3}} \RR^{n_{rb}k}\xrightarrow{\text{reshape}} \RR^{n_{rb}\times k}$$
\begin{itemize}
    \item $L^{rb,d}_1$: $\,\, V_{t,h} \xrightarrow{r=0.3} V_{t,H}$,
    \item $L^{rb,d}_2$: $\,\, V_{t,H} \xrightarrow{r=0.3} V_{t,H}$,
    \item $L^{rb,d}_3$: $\,\,V_{t,H} \xrightarrow{\textcolor{white}{r=0.3}} \RR^{n_{rb}k}\xrightarrow{\text{reshape}} \RR^{n_{rb}\times k}$,
\end{itemize}
where $V_{t,H}=X_{H}^1(\Omega)$ is the finite-element space determined from the computational mesh $\mathcal{T}_{H}$ of tetrahedrons, with stepsize $H>h$ (8$\times$8$\times$3 grid). The network $\minn_{rb,d}(\dist; \hyperd)$ generates a matrix in $\RR^{n_{rb}\times k}$. Each column of such matrix is a vector in the space of the reduced coefficients. These $k=10$ vectors are linearly combined through some hyper-parameters (described below) to obtain the reduced coefficients that best fit the loss function.

The same considerations holds true for the neural network taking the inlet function $\inlets$ as input, with the only difference that in the dense layer $L^{rb,d}_3$ no activation is applied, while the dense layer $L^{rb,\eta}_3$ of the MINN $\minn_{rb,\eta}$ is composed with a hyperbolic tangent function. In conclusion:
\begin{equation*}
     \minn_{rb,d}, \minn_{rb,\eta} : V_{t,h} \xrightarrow{r=0.3} V_{t,H} \xrightarrow{r=0.3} V_{t,H} \xrightarrow{\textcolor{white}{r=0.6}} \RR^{n_{rb}\times 10},
\end{equation*}
giving back two matrices $\mathbb{X}_{\dist} \in \RR^{n_{rb} \times 10}$ and $\mathbb{X}_{\inlets} \in \RR^{n_{rb} \times 10}$, respectively.

For what concerns the physical parameters $\parma \in \RR^3$, a shallow neural network $\ann_{rb,\macro}: \RR^{3} \rightarrow \RR^{10}$ with two dense layers $L^{rb,\micro}_1$ and $L^{rb,\micro}_2$ is introduced: 
%\begin{itemize}
%    \item $L^{rb,\micro}_1 : \RR^3 \xrightarrow{\textcolor{white}{r=0.6}} \RR^{15}$,
%    \item $L^{rb,\micro}_2 : \RR^{15} \xrightarrow{\textcolor{white}{r=0.6}} \RR^{10}$,
%\end{itemize}
$$L^{rb,\micro}_1 : \RR^3 \xrightarrow{\textcolor{white}{r=0.6}} \RR^{15},\quad\quad L^{rb,\micro}_2 : \RR^{15} \xrightarrow{\textcolor{white}{r=0.6}} \RR^{10},$$
where the input layer is composed with a hyperbolic tangent function and the output one has no activation. This leads to the architecture
\begin{equation*}
     \ann_{rb,\macro} : \RR^3 \xrightarrow{\textcolor{white}{r=0.6}} \RR^{15} \xrightarrow{\textcolor{white}{r=0.6}} \RR^{10},
\end{equation*}
that provides as output a vector $\mathbf{y}_{\parma}=\ann_{rb,\macro}(\parma) \in \RR^{10}$.

In order to obtain the approximation of the reduced coefficients, we hinge on the aforementioned data structures and we perform the following linear combination
\[\muurb = (\mathbb{X}_{\inlets} \odot \mathbb{X}_{\dist})\mathbf{y}_{\parma}= \mathbb{X}_{\parmi}\mathbf{y}_{\parma}%=\big(\sum_{j=1}^k \mathbb{X}(:,j)\mathbf{y}_j\big)
.\]
The network $\minn_{rb}$ of equation \eqref{eq:POD_MINN_micro} is trained for at most $200$ epochs, minimizing the loss function
\begin{equation*}
    \mathcal{E}(\minn_{rb};\hyper_{rb}) = \frac{1}{N_{train}} \sum_{i=1}^{N_{train}} \| \VV^T\SS(:,i) - \minn_{rb}(\para^{(i)};\hyper_{rb})\|_{2,\nrb},
\end{equation*}
relying on the L-BFGS optimizer with learning rate equal to 1 and no batching  and applying the same early stopping criterion presented in Section 4.3.
We initialize the corresponding weights and biases with the normalized initialization proposed in \cite{GlorotBengio}. 

\subsubsection{Architectures and training of the closure model (POD-MINN+)}
The augmented POD-MINN approach is carried out exploiting a closure model $\minn_c (\parmi; \hyper_c): \RR^{N_h} \to \RR^{N_h}$ that is fed only with the geometrical parameters $\dist$ and $\inlets$, encoding the microscale features of the microcirculation problem. As in the previous case, we introduce two distinct MINNs that handle each of two inputs, i.e. $\minn_{c,d}$ and $\minn_{c,\eta}$, separating their individual effects on the correction of the POD-MINN model: 
\[\minn_c (\parmi; \hyper_c) = \minn_{c,\eta}(\inlets; \hyper_{c,\inlets}) \odot \minn_{c,d}(\dist; \hyper_{c,d}). \]
Each of these MINNs is built with a single mesh-informed layer:
\begin{itemize}
    \item $L_1^{c,d}$: $V_{t,h} \xrightarrow{r=0.5} V_{t,h}$, activation function $\rho_d(x)=0.4\,\,\tanh(x)$;
    \item $L_1^{c,\eta}$: $V_{t,h} \xrightarrow{r=0.5} V_{t,h}$, activation function $\rho_{\eta}(x)=0.05\,\, \tanh(x)$.
\end{itemize}
The outputs are the vectors $\mathbf{z}_{\dist}=\minn_{c,d} (\dist) \in \RR^{N_h}$ and $\mathbf{z}_{\inlets}=\minn_{c,\eta} (\inlets) \in \RR^{N_h}$, providing the following POD-MINN+ approximation:
\begin{equation*}
  \widetilde{\mathbf{u}}_{\para,h}^{rb} 
  =\VV \muurb + \nuuu
  =\VV\mathbb{X}_{\parmi}\mathbf{y}_{\parma} + (\mathbf{z}_{\inlets} \odot \mathbf{z}_{\dist}) \approx \uuph.
\end{equation*}
To train the closure model the following loss function is minimized:
\begin{align*}
    \mathcal{E}(\minn_c;\hyper_c) = \frac{1}{N_{train}} \sum_{i=1}^{N_{train}} \Bigg(&\| \big(\SS(:,i) - \VV\minn_{rb}(\para^{(i)};\hyper_{rb})\big)
    - \minn_c(\parmi^{(i)};\hyper_c) \|_{2,N_h} + \\ 
    &\| \big(\SS(:,i) - \VV\minn_{rb}(\para^{(i)};\hyper_{rb})\big)
    - \minn_c(\parmi^{(i)};\hyper_c) \|_{\infty,N_h} \Bigg),
\end{align*}
where we have introduced a regularization term to adjust the loss to control the reconstruction error of the solution in the $L_{\infty}$ norm. 

The training consists of at most $200$ epochs with the same early stopping criterion previously mentioned. Also in this case we rely on the L-BFGS optimizer (learning rate equal to 1 and no batching). Thanks to the fact that the closure model builds upon the POD-MINN method that is already trained, the weights and biases collected into the vector $\hyper_c$ are initialized to zero.

It is crucial to stress that for this real-life problem with complex input data, an approach based on training directly a parameter-to-solution neural network (based on MINNs) is not feasible. The combination of the linear projection achieved by the POD-MINN followed by the closure model is the key feature for the success of the method. We need firstly to perform a POD projection into a reduced basis space, so that the closure model acts as a correction of an already trained model. 

\subsubsection{Computational performance}
The speed up provided by the POD-MINN (with or without closure) is significant, unlocking approaches that would be unfeasible for the full order model and allowing to employ the POD-MINN+ model as a surrogate for the numerical solver of the microcirculation problem. The trained parameter-to-solution POD-MINN+ model is able to compute a solution in approximately $0.001$ seconds, while a single run of the FOM requires a higher wall computational time, that varies from $15$ to $35$ minutes, depending on the density of the embedded 1D microstructure.

The gain with respect to the FOM is relevant also if we include the training times.
In average, the training time of the POD/MINN method oscillates between 50 and 100 seconds, depending on the number of basis functions (but the trend is not monotone). The training time for the closure model is higher, although the optimization process consists of a correction of a trained model, and it is almost independent of the number of the POD basis functions, varying from 244 to 276 seconds.

% \begin{figure}[htb]
% \centering
% \includegraphics[height=0.21\textheight]{Images/trainingtimes.png}
% \caption{POD-MINN and closure model training times, plotted with respect to the number of POD basis functions.}
% \label{fig:trainingtimesmicro}
% \end{figure}

% A current bottleneck of this procedure is the wall computational time needed to sample the geometrical parameters in the pre-processing phase. In particular the generation of the vascular networks is a costly process, based on a trial-error algorithm that determines a microstructure with geometrical features that match the requested value of the markers explained in Table \ref{tab:parameters_geom}.

\begin{figure}[htb]
\centering
\includegraphics[width=0.49\textwidth]{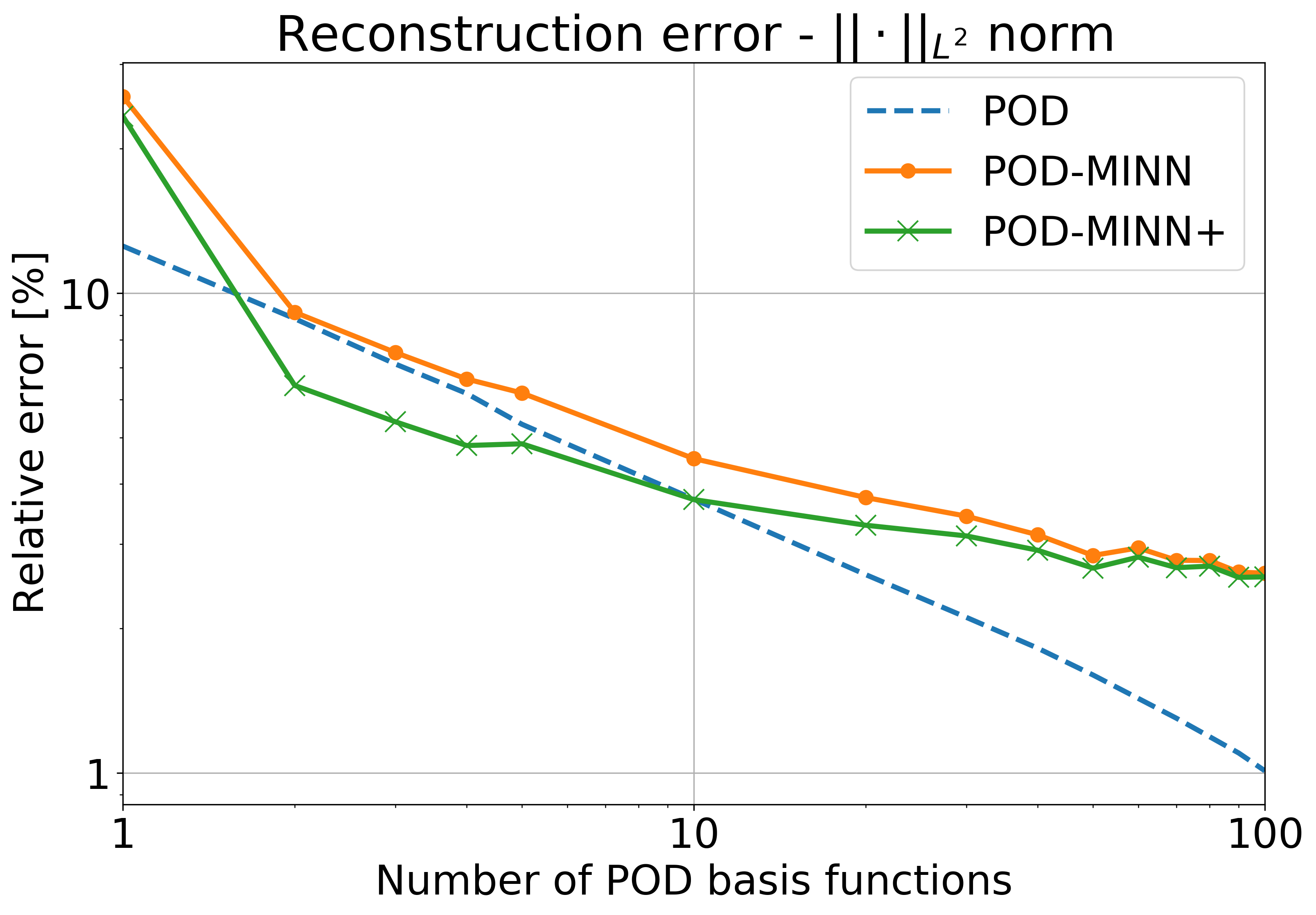}
\includegraphics[width=0.49\textwidth]{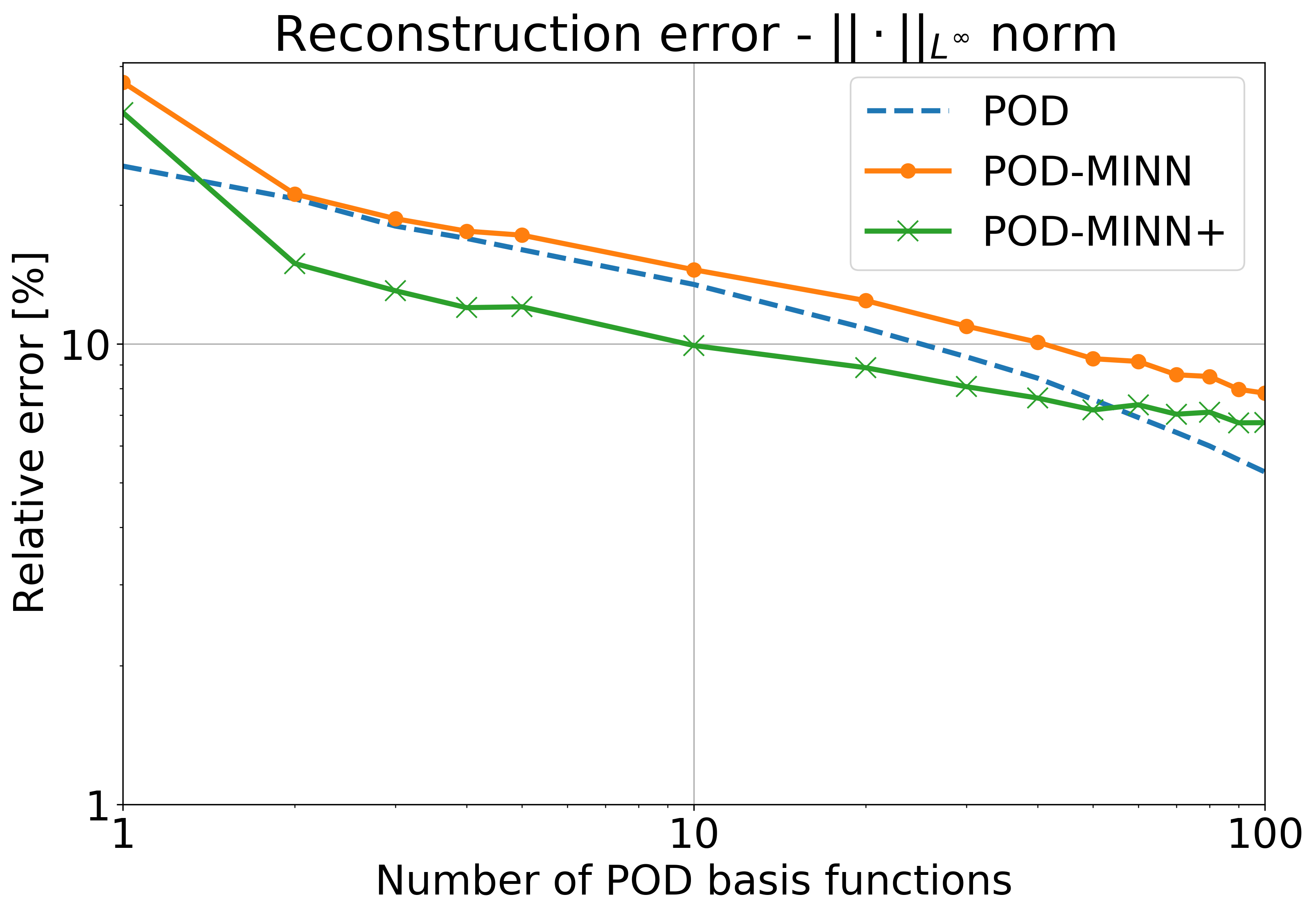}
\caption{These plots compare the error $E_{\tiny POD-MINN}$ and $E_{\tiny POD-MINN+}$ defined in \eqref{eq:EPODMINN}, to the one of the POD projection defined as $E_{POD}$ in \eqref{eq:EPOD}, measured in the norms $p=2,\infty$ computed for the oxygen transfer problem \eqref{eq:oxy}. The case with $p=2$ is reported on the left, the one with $p=\infty$ is on the right. The plots are shown in logarithmic scale.}
\label{fig:microPODMINN+}
\end{figure}

\subsection{Numerical results: comparing linear and nonlinear model order reduction}

In order to measure the approximation properties of the POD-MINN+ method, we analyze the POD projection error $E_{POD}$, defined as in \eqref{eq:EPOD}, and the reconstruction errors $E_{\tiny PODMINN}$ and $E_{\tiny PODMINN+}$. We compare the behaviour of each error varying the number of POD basis functions $\nrb$ in Figure \ref{fig:microPODMINN+}. 

From a first analysis, we notice the good performance of both methods, being able to capture featuring good approximation results for the $L^2$ norm of the errors with few POD basis functions. Indeed oxygen transfer is diffusion dominated and consequently it is well suited for an approximation methods based on POD. This is confirmed by the fast decay of the singular values of $\mathbb{S}$ described before. 
On the other hand, a significant difference is detected between the $L^2$ and the $L^{\infty}$ norm, as the POD-MINN+ shows overall higher accuracy with respect to POD and POD-MINN. The latter approach fails to retrieve the local effects due to the microscale data. Instead, the closure model is able to exploit the information about the microscale geometry to improve the smallest scales, provided that they are resolved by the FOM. 

The efficacy of the POD/MINN method is more significant with low number of POD basis functions, as the closure model approximates the high-frequency modes that are neglected at POD level: see Table \ref{tab:errors_micro_5_10_20_bases} to gain more comprehensive quantitative insights about each proposed method.

\begin{table}[h!t]
    \centering
    \renewcommand{\arraystretch}{1.75} % Adjust the vertical spacing between rows
    \setlength{\tabcolsep}{7pt} % Adjust the horizontal spacing between columns
    \begin{tabular}{||l|c c|c c|c c||}
        \cline{2-7}
        \multicolumn{1}{c|}{} & \multicolumn{6}{c|}{\textsc{Test Errors}} \\ [1ex]
        \cline{2-7}
        \multicolumn{1}{c|}{} & \multicolumn{2}{c|}{$5$ \textbf{POD Modes}} & \multicolumn{2}{c|}{$10$ \textbf{POD Modes}} & \multicolumn{2}{c|}{$20$ \textbf{POD Modes}} \\
        \hline
        \textsc{Method} & \textsc{$\|\cdot\|_{2,N_h}$} & \textsc{$\|\cdot\|_{\infty,N_h}$} & \textsc{$\|\cdot\|_{2,N_h}$} & \textsc{$\|\cdot\|_{\infty,N_h}$} & \textsc{$\|\cdot\|_{2,N_h}$} & \textsc{$\|\cdot\|_{\infty,N_h}$} \\ \hline
        \textbf{POD} & $5.33\%$ & $16.00\%$ & $3.72\%$ & $13.45\%$ & $2.59\%$ & $10.80\%$ \\ [1ex]
        \textbf{POD-MINN} & $6.19\%$ & $17.21\%$ & $4.52\%$ & $14.48\%$ & $3.75\%$ & $12.41\%$ \\ [1ex]
        \textbf{POD-MINN}+ & $4.85\%$ & $12.03\%$ & $3.71\%$ & $9.91\%$ & $3.28\%$ & $8.87\%$ \\ [1ex] \hline
    \end{tabular}
    \caption{Comparison of POD, POD-MINN, and POD-MINN+ with respect to different choices of the number of reduced basis functions. We report the errors $\|\uuph-\muuphrb\|_{p,N_h}$ and $\|\uuph-\nuuphrb\|_{p,N_h}$, for the POD-MINN and POD-MINN+ respectively, in the Euclidean ($p=2$) and maximum norm ($p=\infty$).}
    \label{tab:errors_micro_5_10_20_bases}
\end{table}

\begin{figure}[h!tb]
\centering
\includegraphics[height=0.25\textheight]{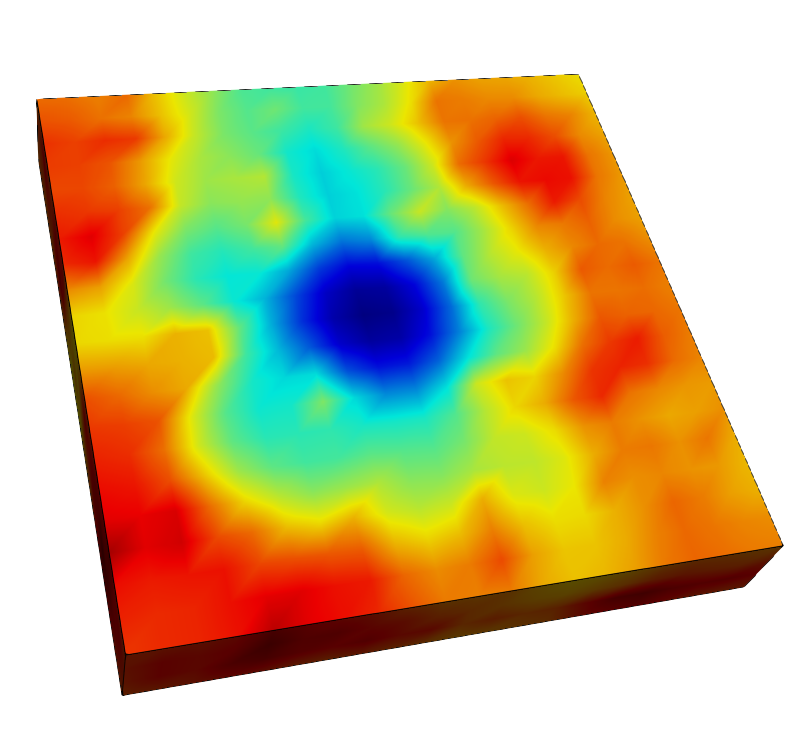}
\includegraphics[height=0.25\textheight]{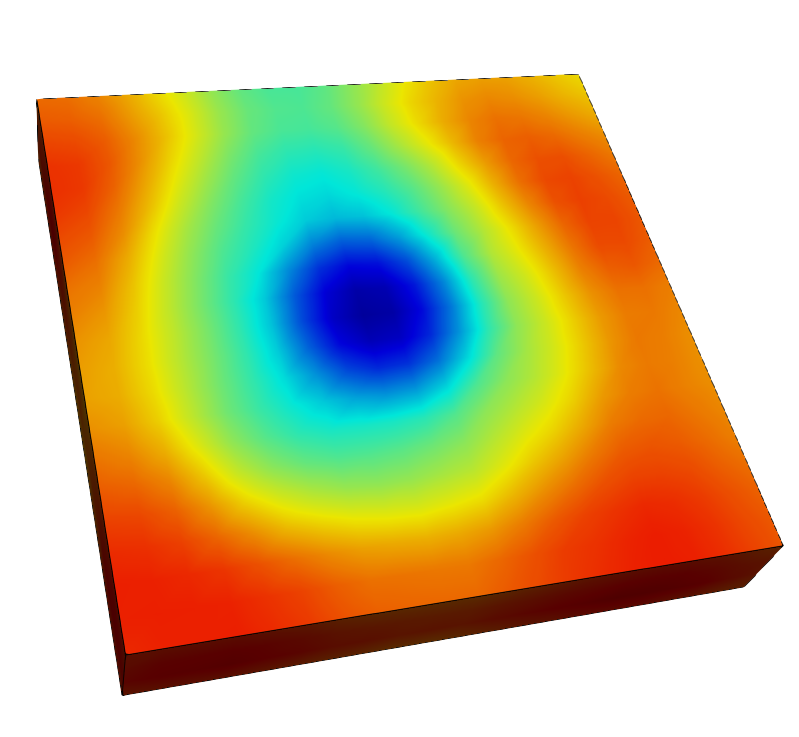} \\
a) FOM solution \textcolor{white}{onstruction} \hspace{1cm} b) POD projection\textcolor{white}{nstruction} \\
\includegraphics[height=0.25\textheight]{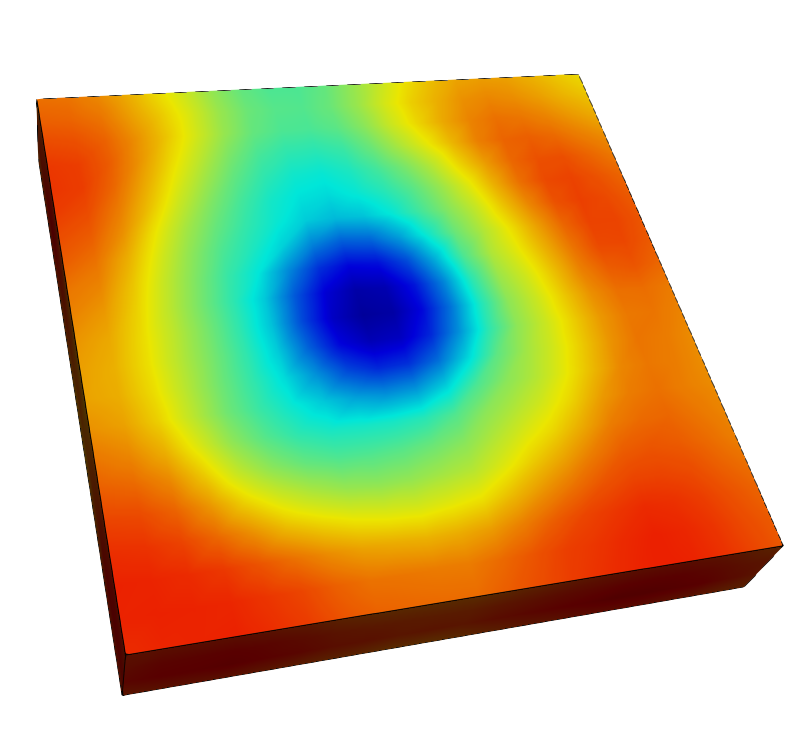}
\includegraphics[height=0.25\textheight]{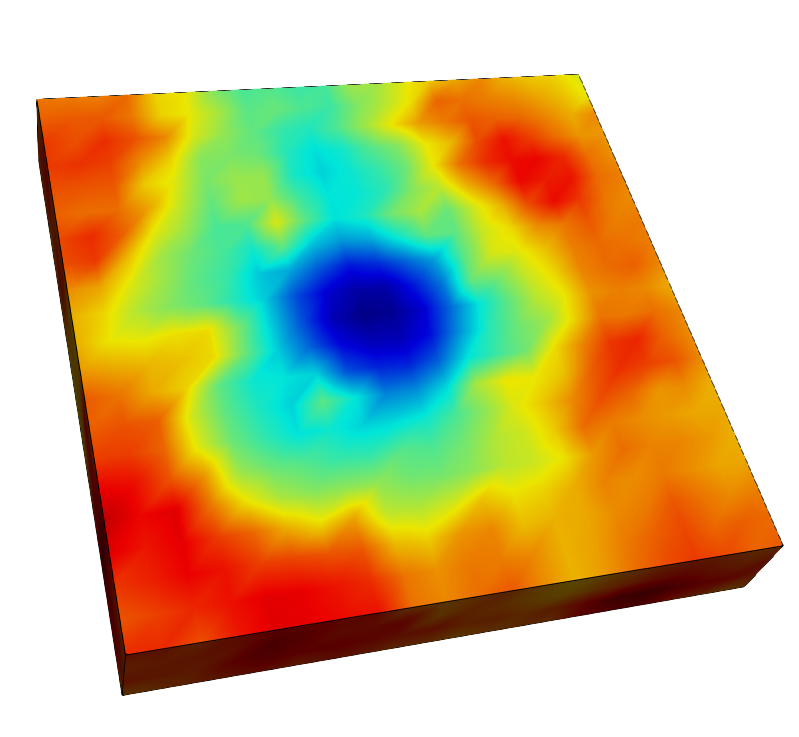} \\
c) POD-MINN approximation \hspace{1cm} d) POD-MINN+ approximation \\
\includegraphics[height=0.07\textheight]{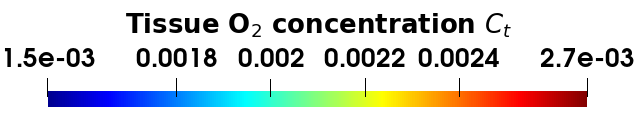} \\
\caption{Plots of the high-fidelity FOM solution in a) is compared with the corresponding POD projection in b) and the reconstructions with the POD-MINN and POD-MINN+ approaches in c) and d) respectively: through the former method, given a small number of POD basis functions ($n_{rb} = 10$), it is possible to retrieve only the global features of the solution at the macroscale, while the local effects associated to the macroscale are captured thanks to the closure model.}
\label{fig:reconstruction_comparison}
\end{figure}

This analysis becomes more evident when we compare in Figure \ref{fig:reconstruction_comparison} the visualization of a FOM solution with the three reconstructions, obtained respectively using POD, POD-MINN and POD-MINN+. It is clear that the POD and POD-MINN methods are over-diffusive, failing in representing the details of the solution at the smallest scale that are filtered out due to the projection on the reduced basis space of only $n_{rb}=10$ basis functions. On the other hand, the POD-MINN+ ensures the microstructures are correctly captured, and yields a much more appreciable description of the true behaviour of the solution over the entire domain.

\section{Conclusions}\label{sec:conclusions}

Reduced order modeling plays a crucial role in approximating the behaviour of continuous media with microstructure. These complex systems exhibit intricate geometries and intricate physics, making their direct numerical simulation computationally expensive and time-consuming. 

Reduced order modeling offers an efficient approach to capture the essential features of the system while significantly reducing computational costs. However, well established techniques such as proper orthogonal decomposition (POD) and reduced basis methods, exploiting the projection of functions belonging to high-dimensional spaces onto a low-dimensional subspace, fail in preserving the typical high frequency modes of the true solution in presence of a microstructure. On the other hand, neural networks possess a remarkable ability to approximate functions in high dimensions. This flexibility allows neural networks to capture intricate patterns and relationships even in high-dimensional spaces. Leveraging on these properties, we have developed a correction of the POD-Galerkin approximation that restores the fine scale details into the reduced solution.

We perform this task in two steps: first we use MINNs to approximate the map from the parameters of the problem to the reduced POD coefficients, yielding the POD-MINN method; second we enhance the approach by including an additive closure model, that is a correction term ultimately providing the POD-MINN+ method. The whole procedure can be defined and successfully trained thanks to a new family of neural network architectures, called mesh-informed neural networks (MINNs), which take advantage of the mesh structure by considering connectivity information between neighboring points or elements in the mesh. This information can be used to design specialized layers that exploit the spatial relationships within the mesh. By incorporating the mesh structure into the network architecture, the model can learn more effectively from the data and capture spatial dependencies that would be difficult to capture using traditional neural network architectures. 

The resulting ROM is accurate, efficient, and non intrusive, thus applicable to many scenarios where a FOM capable to capture the effect of microstructures is available but extremely expensive to query. The POD-MINN+ strategy not only accelerates simulations, but is also potentially able to enhance parametric studies such as sensitivity analysis and uncertainty quantification, making it a valuable tool for understanding the role of microstructure in many physical systems.

\section*{Acknowledgments}
%The present research is part of the activities of “Dipartimento di Eccellenza 2023-2027" 
The present research is part of the activities of project Dipartimento di Eccellenza 2023-2027, funded by MUR, and of project FAIR (Future Artificial Intelligence Research) project, funded by the NextGenerationEU program within the PNRR-PE-AI scheme (M4C2, Investment 1.3, Line on Artificial Intelligence). The authors are members of Gruppo Nazionale per il Calcolo Scientifico (GNCS) of Istituto Nazionale di Alta Matematica (INdAM).



\bibliographystyle{elsarticle-num} 
\bibliography{bibliography.bib}

\end{document}